%% file: log.tex
\documentclass[12pt]{article}

\usepackage{graphicx}
\usepackage{latexsym,amsmath,amsfonts,amscd, amsthm}
\usepackage{color}
\usepackage{bm}
\usepackage[titletoc]{appendix}
\usepackage{multirow}
\numberwithin{equation}{section}

\newtheoremstyle{plainNoItalics}{}{}{\normalfont}{}{\bfseries}{.}{ }{}

\theoremstyle{plain}
\newtheorem{thm}{Theorem}[section]

\theoremstyle{plainNoItalics}

\newtheorem{defn}[thm]{Definition}
\newtheorem{rem}[thm]{Remark}
\newtheorem{prop}[thm]{Proposition}
\newtheorem{exa}[thm]{Example}

\newcommand{\f}{\frac}

\newcommand{\beq}{\begin{equation}}
\newcommand{\eeq}{\end{equation}}
\newcommand{\beqa}{\begin{eqnarray}}
\newcommand{\eeqa}{\end{eqnarray}}
\newcommand{\bit}{\begin{itemize}}
\newcommand{\eit}{\end{itemize}}
\newcommand{\bedef}{\begin{defn}}
\newcommand{\edefn}{\end{defn}}
\newcommand{\bpro}{\begin{prop}}
\newcommand{\epro}{\end{prop}}

\newcommand{\Dx}{\Delta x}

\newcommand{\pad}[2]{\frac{\partial #1}{\partial #2}}

 \newcommand{\rf}[1]{(\ref{#1})}

 \newcommand{\sumls}{\sum_{\ell=1}^s}
 \newcommand{\suml}[1]{\sum_{\ell=1}^{#1}}
 \newcommand{\tC}{\tilde{C}}
 \newcommand{\tc}{\tilde{c}}

\def\Box{\mbox{ }\rule[0pt]{1.5ex}{1.5ex}}

\begin{document}

\baselineskip=1.8pc


\input{title}

\newpage

\input{intro}

\input{cons}

\input{ode_stab}

\input{pde_stab}
\input{numerics}
\input{conclusion}

\bibliographystyle{siam}
\bibliography{refer}

\end{document}

%% file: title.tex
\begin{center}
{\bf
Conservative Multi-Dimensional Semi-Lagrangian Finite Difference Scheme: Stability and Applications to the Kinetic and Fluid Simulations
}
\end{center}

\vspace{.2in}
\centerline{
Tao Xiong  \footnote{School of Mathematical Sciences, Fujian Provincial Key Laboratory of Mathematical Modeling and High-Performance Scientific Computing, Xiamen University, Xiamen, Fujian, P.R. China, 361005. Email: txiong@xmu.edu.cn}
Giovanni Russo \footnote{Department of Mathematics and Computer Science, University of Catania, Catania, 95125, E-mail: russo@dmi.unict.it} 
Jing-Mei Qiu\footnote{Department of Mathematics, University of Houston,
Houston, 77004. E-mail: jingqiu@math.uh.edu. Research supported by
Air Force Office of Scientific Computing grant FA9550-16-1-0179, NSF grant DMS-1522777 and University of Houston.}
}

\bigskip
\centerline{\bf Abstract}
\bigskip
In this paper, we propose a mass conservative semi-Lagrangian finite difference scheme for multi-dimensional problems {\em without dimensional splitting}.
The semi-Lagrangian scheme, based on tracing characteristics backward in time from grid points, does not necessarily conserve the total mass.
To ensure mass conservation, we propose a conservative correction procedure based on a flux difference form. Such procedure guarantees local mass conservation, while introducing time step constraints for stability. We theoretically investigate such stability constraints from an ODE point of view by assuming exact evaluation of spatial differential operators and from the Fourier analysis for linear PDEs.

The scheme is tested by classical two dimensional linear passive-transport problems, such as linear advection, rotation and swirling deformation. The scheme is applied to solve the nonlinear Vlasov-Poisson system using a a high order tracing mechanism proposed in [Qiu and Russo, 2016]. Such high order characteristics tracing scheme is generalized to the nonlinear guiding center Vlasov model and incompressible Euler system. The effectiveness of the proposed conservative semi-Lagrangian scheme is demonstrated numerically by our extensive numerical tests.
\vfill


\newpage

%% file: intro.tex
\section{Introduction}
\label{sec1}
\setcounter{equation}{0}
\setcounter{figure}{0}
\setcounter{table}{0}

Semi-Lagrangian (SL) schemes have been used extensively in many areas of science and engineering, including weather forecasting \cite{staniforth1991semi, lin1996multidimensional, guo2014conservative}, kinetic simulations \cite{filbet2001conservative, guo2013hybrid} and fluid simulations \cite{pironneau1982transport, xiu2001semi}, interface tracing \cite{enright2005fast, strain1999semi}, etc. The schemes are designed to combine the advantages of Eulerian and Lagrangian approaches. In particular, the schemes are build upon a fixed set of computational mesh. Similar to the Eulerian approach, high spatial resolution can be realized by using high order interpolation/reconstruction procedures or by using piecewise polynomial solution spaces. On the other hand, in each time step evolution, the scheme is designed by propagating information along characteristics, relieving the CFL condition. Typically, the numerical time step size allowed for an SL scheme is larger than that of an Eulerian approach, leading to gains in computational efficiency.

Among high order SL schemes, depending on solution spaces, different classes of methods can be designed. For example, a finite difference scheme evolves point-wise values and realizes high spatial resolution by high order interpolation procedures \cite{xiu2001semi, qiu2010conservative}, a finite volume scheme considers integrated cell-averages with high order reconstruction procedures \cite{lin1996multidimensional, crouseilles2010conservative}, while a finite element method has piecewise continuous or discontinuous polynomial functions as its solution space \cite{pironneau1982transport, morton1988stability, qiu2011positivity, rossmanith2011positivity, guo2014conservative}. Each class of the above mentioned SL methods has its own advantages. For example, the finite element method is more flexible with the geometry and handling boundary conditions, while the finite difference and finite volume schemes could perform better in resolving solution structures with sharp gradients, e.g. by using a weighted essentially non-oscillatory (WENO) procedure. To compare finite difference and finite volume schemes, finite volume scheme is often considered more physically relevant and the local mass conservation can be built up in a natural way; while the finite difference scheme is more flexible and computationally efficient for high-dimensional problems, if one consider schemes of third order or higher.

In this paper, we consider SL finite difference scheme with {\em local mass conservation property}. In fact, many existing SL
finite difference schemes are built based on tracing characteristics backward in time together with a high order interpolation procedure \cite{carrillo2007nonoscillatory}. Typically such schemes do not have local mass conservation property, which is fine for some certain applications. However, for applications in weather forecasting or in kinetic simulations, ignoring local mass conservation could lead to significant loss of total mass, especially when the solution with sharp gradients becomes under-resolved by the computational mesh \cite{huot2003instability}.

There have been many attempts to preserve the mass conservation of an SL finite difference scheme with large time stepping sizes, e.g. \cite{qiu2010conservative, qiu2011conservative_jcp}. However, they are mostly designed for 1D problems by taking advantage of some special features in a 1D setting. Their generalization to high dimensional problems often relies on dimensional splitting which is subject to splitting errors. In this paper, we propose and investigate a truly multi-dimensional approach without dimensional splitting errors. To build an SL finite difference scheme with local mass conservation, one essential framework that we propose to work with is the flux-difference form. However, by working with the flux difference form, one often observe time step constraint for numerical stability. Note that unlike the Eulerian approach, such time step constraint does not come from the CFL condition (i.e. numerical domain of dependence should include the physical domain of dependence), but from numerical stability one employed for temporal integration.
As far as we are aware of, there  is little work in quantifying the stability constraint and in optimizing the numerical strategies balancing stability, accuracy and computational efficiency. This paper aims to fill such gap, by understanding such time step constraints.
In particular, we investigate the stability of time integration schemes based on a linear stability analysis around the imaginary axis in the complex plane, assuming the spatial differentiation is exact. We optimize quadrature rules for time integration by maximizing the stability interval along the imaginary axis. We further employ Fourier analysis to study the numerical stability of a fully-discretized scheme. The schemes are applied to 2D passive-transport problems, as well as to nonlinear Vlasov-Poisson (VP) by using a high order characteristics tracing scheme proposed in \cite{qiu_russo_2016}. Further more, we apply the scheme to the nonlinear guiding center Vlasov system and the incompressible Euler system in vorticity stream function formulation, for which we propose a high order characteristics tracing scheme following the idea in \cite{qiu_russo_2016}.
Finally, we would like to mention a few of our previous work related to stability of SL finite difference schemes in flux-difference form. In \cite{qiu2011conservative_jcp}, a special treatment is introduced to relieve the time step constraint for 1D passive transport problems. However, such treatment is not possible for general high dimensional problems. In \cite{christlieb2014high}, the time step constraint is studied by Fourier analysis for an SL finite difference scheme coupled with integral deferred correction framework.

The paper is organized as follows. The SL finite difference scheme in flux-difference form is described in Section \ref{sec: cons}. The stability of time integration with quadrature rule is investigated in Section \ref{ssec: temp}, assuming {\em exact} evaluation of spatial differentiation operators.
We also optimize temporal integration rules. In Section \ref{ssec: spat}, we study the numerical stability of a fully discretized scheme by Fourier analysis. In Section \ref{sec5}, numerical tests are performed for 2D linear passive-transport problems. In Section \ref{sec6}, we apply the scheme to the nonlinear VP system, the guiding center Vlasov equation and incompressible Euler system.

%% file: cons.tex
\section{A mass conservative SL finite difference scheme}
\label{sec: cons}

In this section, we describe an SL finite difference scheme based on a flux-difference form to locally preserve mass. The scheme starts from a standard non-conservative procedure with backward characteristics tracing and high order spatial interpolation. Then a conservative correction is performed by a flux-difference formulation. We describe the scheme in a 1D linear setting, noting that its extension to nonlinear and high dimensional problems is straightforward, as long as characteristics can be properly traced backward in time, e.g. see our numerical examples in Section \ref{sec5}.

We consider a 1D linear advection equation,
\begin{equation}
\pad{f}{t} + \pad{f}{x} = 0, \quad f(x,0) = f^0(x), \quad x\in [-\pi, \pi].
\label{eq:scalar2}
\end{equation}
For simplicity, we assume a periodic boundary condition.
We assume a uniform discretization in space with $x_j = j\Dx$, $j=1\ldots,n_x$ and  let $f^n_j$ be an approximation of the solution at time $t^n$ and position $x_j$. We describe below the conservative SL procedure to update $\{f^{n+1}_j\}_{j=1}^{n_x}$ from $\{f^n_j\}_{j=1}^{n_x}$.

In an Eulerian finite difference method, typically one would firstly approximate the spatial derivative by a flux difference form to ensure mass conservation, then the system of ODEs will be evolved in time by a high order numerical integrator such as the Runge-Kutta (RK) method via the method of lines. In the SL setting, however, we propose to perform the time integration based on quadrature rules first,
\beq
\label{eq: conser}
f^{n+1}_j = f^n_j - \frac{\partial}{\partial x}\left(\int_{t^n}^{t^{n+1}} f(x, t)dt\right)|_{x=x_j} \approx
f^n_j - \frac{\partial}{\partial x} (\mathcal{F}(x))|_{x_j},
\eeq
where we let
\beq
\label{eq: quadr}
\mathcal{F}(x) \doteq  \sum_{\ell=1}^{s} f(x, t^n+c_\ell \Delta t) b_\ell \Delta t
\eeq
as a quadrature approximation of $\int_{t^n}^{t^{n+1}} f(x, t)dt$. Here $(c_\ell,b_\ell)$, $\ell = 1,\ldots,s$ are the nodes and weights of an accurate quadrature formula and $f(x_j, t^n+c_\ell \Delta t)$, $\ell=1\cdots s$ (we call it stage values) can be approximated via a non-conservative SL scheme via backward characteristics tracing and high order spatial interpolation. For the linear equation \eqref{eq: conser}, $f(x_j, t^n+c_\ell \Delta t)$ can be traced back along characteristics to $t^n$ at $f(x_j-c_\ell \Delta t/\Delta x, t^n)$, whose value can be obtained via a WENO interpolation from neighboring grid point values  $\{f^n_j\}_{j=1}^{n_x}$, see our description for different interpolation procedures in Section~\ref{ssec: spat}.
Then a conservative scheme, based on a flux-difference form, can be proposed in the spirit of the work by Shu and Osher \cite{shu1989efficient}. In particular, the scheme can be formulated as
\beq
\label{eq: cons_update}
f^{n+1}_j = f^n_j - \frac{1}{\Delta x} (\hat{F}_{j+\frac12} - \hat{F}_{j-\frac12}),
\eeq
where $\hat{F}_{j+\frac12}$ comes from WENO reconstruction of fluxes from $\{\mathcal{F}_j\}_{j=1}^{nx}$ with $\mathcal{F}_j \doteq \mathcal{F}(x_j)$.
We refer to \cite{shu2009high} for the basic principle and detailed procedures of WENO reconstruction. Also, Section \ref{ssec: spat} provides detailed discussions on different reconstruction procedures. It can be shown that the mass conservation is locally preserved due to the flux difference form \eqref{eq: cons_update}.

Such conservative correction procedure can be directly generalized to problems with non-constant velocity fields in a multi-dimensional setting without any difficulty, e.g. rotation and swirling deformation. In additional to the procedures described above, a high order ODE integrator such as a Runge-Kutta method can be employed to locate the foot of a characteristic accurately.  For example, we consider a 2D problem with a prescribed velocity field $a(x, y, t)$ and $b(x, y, t)$
\[
f_t + \left(a(x, y, t) f\right)_x + \left(b(x, y, t) f\right)_y = 0.
\]
Let the set of grid points
\beq
\label{eq: 2dgrid}
x_1 < \cdots < x_i < \cdots < x_{n_x}, \quad y_1 < \cdots < y_j < \cdots < y_{n_y}
\eeq
be a uniform discretization of a 2D rectangular domain with $x_i = i \Delta x$ and $y_j = j \Delta y$.
The foot of characteristic emanating from a 2D grid point, say $(x_i, y_j)$ at $t^{\ell}\doteq t^n+c_\ell \Delta t$ can be located by solving the following final-value problem accurately with a high order Runge-Kutta method,
\beq
\label{eq: char1}
\frac{dx}{dt} = a(x, y, t), \quad \frac{dy}{dt} = b(x, y, t), \quad x(t^\ell) = x_i, \quad y(t^\ell) = y_j.
\eeq
Once the foot of characteristic located, say at $(x^\star_i, y^\star_j)$, then $f(x_i, y_j, t^{\ell})$ can be evaluated by approximating $f(x^\star_i, y^\star_j, t^n)$ via a high order 2D WENO interpolation procedure \cite{shu2009high}. A 2D conservative scheme based on a flux-difference form can be formulated as
\beq
\label{eq: cons_update}
f^{n+1}_{ij} = f^n_{ij} - \frac{1}{\Delta x} (\hat{F}_{i+\frac12, j} - \hat{F}_{i-\frac12, j})- \frac{1}{\Delta y} (\hat{G}_{i, j+\frac12} - \hat{G}_{i, j-\frac12}),
\eeq
where $\hat{F}_{i\pm\frac12, j}$ comes from WENO reconstruction of fluxes from $\{\mathcal{F}_{ij}\}_{i=1}^{n_x}$ for all $j$ with
\[
\mathcal{F}_{ij} \doteq \mathcal{F}(x_i, y_j) \approx \Delta t \sum_{\ell=1}^{s} f(x_i, y_j, t^n+c_\ell \Delta t) b_\ell.
\]
The procedure for WENO reconstruction is the same as the 1D case for all $j$ and we again refer to the review paper \cite{shu2009high}. Similarly, $\hat{G}_{i, j\pm\frac12}$ comes from WENO reconstruction of fluxes from $\{\mathcal{F}_{ij}\}_{j=1}^{n_y}$ for all $i$.

To generalize the conservative SL scheme to nonlinear systems, a problem-dependent high order characteristics tracing procedure needs to be designed for solving the final-value problem in the form of equation \eqref{eq: char1}, but with the velocity field depending on the unknown function $f$. In many cases, a high order Runge-Kutta method could not be directly applied. In \cite{qiu_russo_2016}, a high order multi-dimensional characteristics tracing scheme for the VP system is proposed and can be applied in the above proposed conservative SL framework.  In Section~\ref{sec6} we present numerical results and generalize the characteristics tracing procedure for the VP system to a guiding center Vlasov system and incompressible Euler system in vorticity stream function formulation.

We close this section by making the following remark to motivate our discussions in the following two sections. There are two sources in the scheme formulation that contribute to the stability issue of the above proposed SL scheme. One is the discretization by the quadrature rule \eqref{eq: quadr}. This part of stability is viewed as an ODE stability (assuming exact evaluation of spatial operators) and is investigated carefully in Section \ref{ssec: temp}. The other source can be explained by observing the following situation: if one changes the time stepping size slightly (could be arbitrary small), the root of characteristics $x_j-c_\ell \Delta t/\Delta x$ could come from a different grid cell, leading to a different interpolation stencil in the implementation. This aspect is associated with spatial discretization and is investigated in Section \ref{ssec: spat}.

%% file: ode_stab.tex
\section{Temporal discretization and stability.}
\label{ssec: temp}
\setcounter{equation}{0}
\setcounter{figure}{0}
\setcounter{table}{0}

\subsection{Linear stability functions and stability regions}
We first investigate the linear stability of quadrature rules for temporal discretization \eqref{eq: quadr} in an ODE setting, by assuming an {\em exact} evaluation of spatial derivative in eq.~\eqref{eq: conser}. In particular, we look for the evolution of a Fourier mode, identified by a Fourier variable $\xi \in [-\pi, \pi]$, assuming exact evaluation of spatial interpolation and reconstruction procedure mentioned above.
Such a discrete Fourier mode at time $t^n = n \Delta t$, will be denoted by,
\[
f_\xi^n(x) = (Q(\xi))^n e^{\mathbf{i} x \xi/\Delta x}, \quad \mathbf{i} = \sqrt{-1},
\]
where $Q(\xi)$ is the amplification factor associated with $\xi$. After plugging such ansatz into the scheme with $c_\ell$ and $b_\ell$, $\ell=1, \cdots s$ for temporal discretization, we obtain
\beq
\label{eq: Q}
Q(\xi) = 1- \mathbf{i} \xi \sum_{\ell=1}^s b_\ell e^{-\mathbf{i}c_\ell\xi}.
\eeq
The scheme is stable if
\[
|Q(\xi)|\le 1, \quad \forall \xi\in[-\pi, \pi].
\]
Such stability property is closely related to the linear stability of the quadrature rule, which can be studied by the stability region for a scalar linear ODE,
\[
z' = y z, \quad z(0)=1, \quad \forall y \in \mathbb{C}.
\]
Considering the quadrature rule with $c_\ell$ and $b_\ell$, $\ell=1, \cdots s$, the associated stability function is
\beq
\label{eq: R}
{R}(y) = 1 + y\sum_{\ell=1}^s b_\ell e^{c_\ell y},
\eeq
with which the stability region can be drawn by the set $\{y\in\mathbb{C}; |R(y)|\le1\}$. Comparing equations~\eqref{eq: Q} and ~\eqref{eq: R}, one has $Q=R(-\mathbf{i}\xi)$, with $\xi\in[-\pi,\pi]$. Thus the stability of a quadrature rule in a conservative SL scheme for a linear advection equation is closely related to the stability on the imaginary axis. In order to guarantee stability, we look for the largest interval $I^*\doteq[-y^*,y^*]$ of the
imaginary axis such that $|R(\mathbf{i}y)| \leq 1$, $\forall y\in I^*$. The bound
\beq
y^*/\pi
\label{eq: cfl_ode}
\eeq quantifies the maximum CFL number for the SL scheme that guarantees
stability.

\begin{figure}
\centering
\includegraphics[width=3.0in]{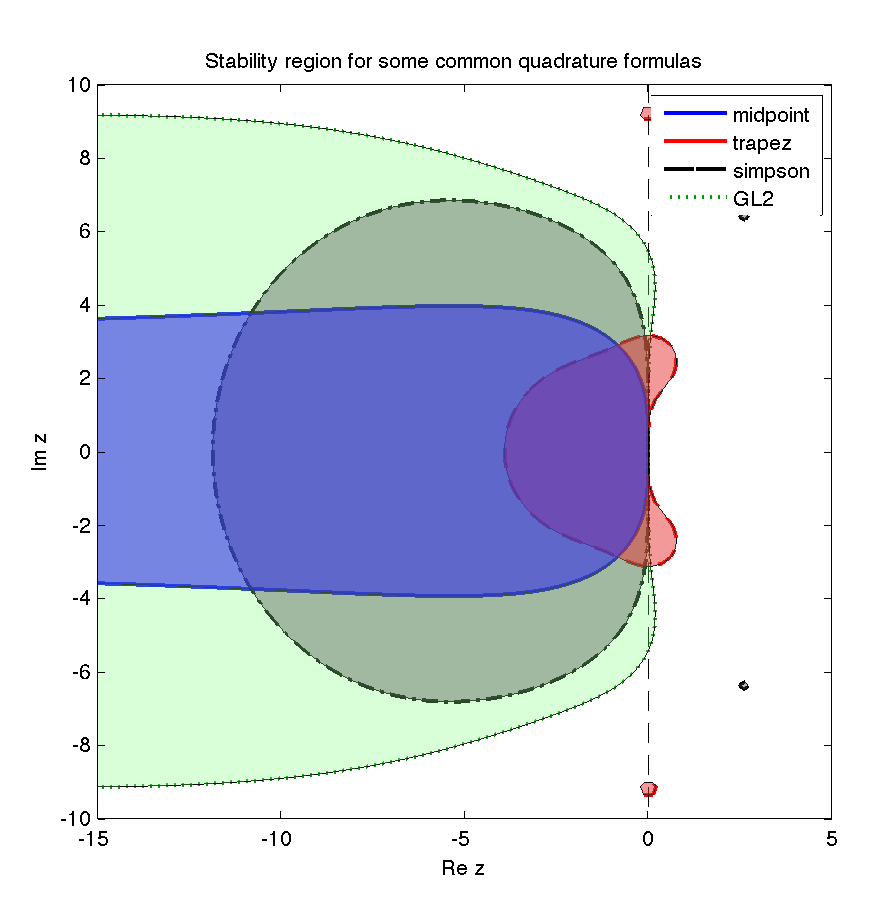}
\includegraphics[width=3.0in]{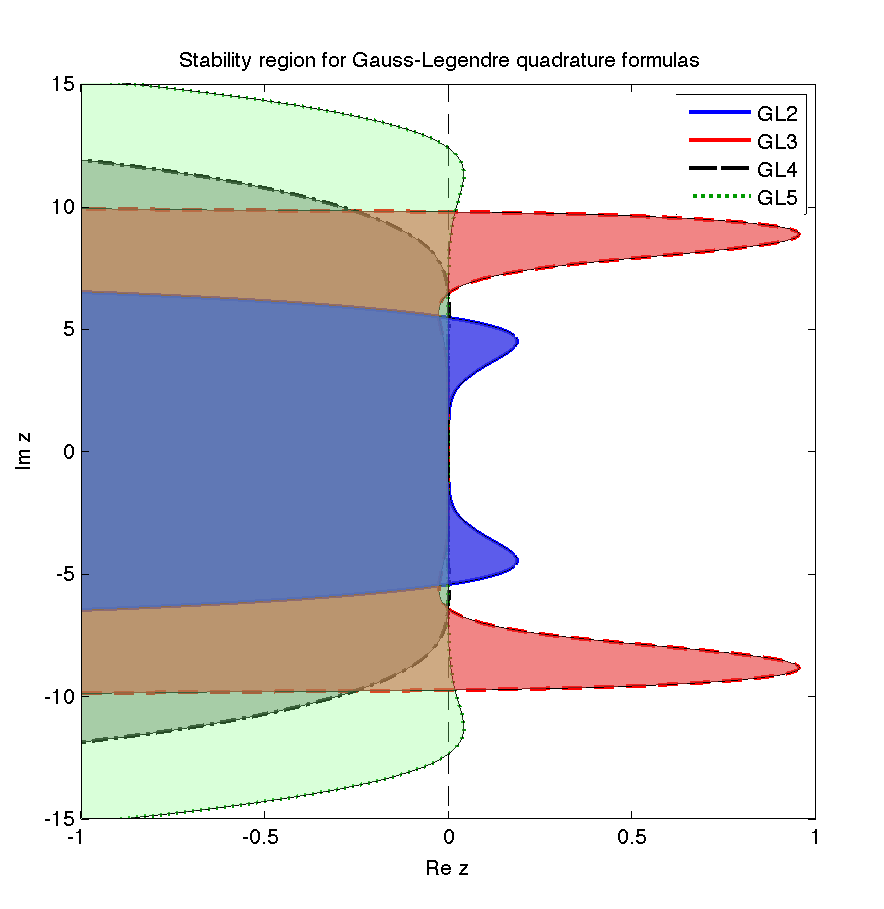}
\caption{Stability region for midpoint, trapezoidal, Simpson's
  and two-point Gauss-Legendre rules (left) and GLs with $s=2, 3, 4, 5$ (right).}
\label{fig:stab_reg}
\end{figure}

\begin{figure}
\centering
\includegraphics[width=3.2in]{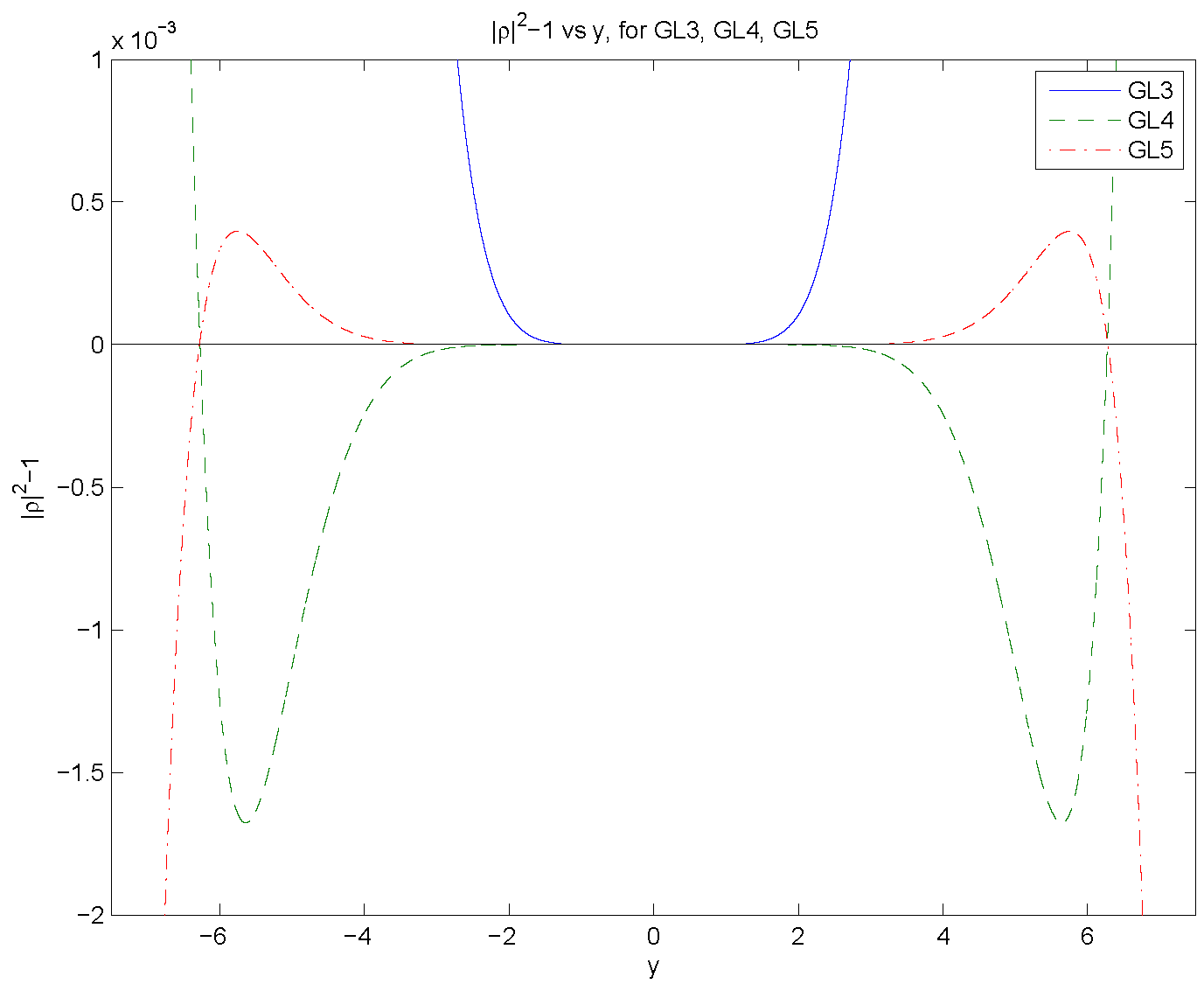}
\caption{$R(\mathbf{i}y)^2-1$ vs $y$ for GL3, GL4, GL5 formulas. GL3 and GL5 is observed to be unstable.}
\label{fig:glstab}
\end{figure}

Below, we report the stability regions for the following commonly used quadrature rules in the left panel of Fig.~\ref{fig:stab_reg}.
\begin{enumerate}
\item midpoint: $q_1 = 1/2$, $w_1=1$.
\item trapezoidal: $q_1 = 0, q_2 = 1$, $w_1 = w_2 = 1/2$.
\item Simpson: $q_1=0, q_2=1/2, q_3=1$, $w_1 = w_3=1/6, w_2=2/3$.
\item two-point Gauss-Legendre formulas (GL2): $q_1=\frac12-\frac{1}{2\sqrt{3}}, q_2=\frac12 +\frac1{2\sqrt{3}}$, $w_1=w_2=1/2$.
\end{enumerate}
As it is apparent from the plot, midpoint and Simpson's rule do not include a portion of the imaginary axis, while the trapezoidal rule and the two-point Gauss-Legendre rule do. The boundary of the stability region of the trapezoidal rule intersects the imaginary axis at $\pi$, and therefore the maximum
CFL number that guarantees linear stability for the conservative scheme is $1$. The two point Gauss-Legendre quadrature formula provides a wider
stability interval, since in this case $y^*\approx 5.43$, giving a maximum CFL number of approximately $5.43/\pi=1.72$.
Higher order Gauss-Legendre quadrature formulas, hereafter denoted by GL$s$, where $s$ indicates the number of nodes, may provide wider stability interval, as is illustrated in the right panel of Fig. \ref{fig:stab_reg}. To better appreciate the stability region, we plot in Fig. \ref{fig:glstab} $\rho^2-1$ as a function of $y$.  GL4 is observed to have better stability property, as $\rho^2-1 \le0$ for an interval with boundary $y^*\approx 6.2765$ leading to a maximum CFL number of approximately $6.2765/\pi=1.99$. Gauss-Legendre rule formulas with odd number of points, such as GL3 and GL5, are unstable near the origin, see the right panel of Fig. \ref{fig:stab_reg} as well as Fig. \ref{fig:glstab}.

\subsection{Maximize the stability interval on imaginary axis}

In order to analyze the stability of quadrature formulas,
let us consider the expression $R(\mathbf{i}y)$ from eq.~\rf{eq: R}, and write it in the form
\begin{equation}
R(\mathbf{i}y) = 1+\mathbf{i}y(C_s(y) + i S_s(y)) = 1-yS_s(y) + \mathbf{i}y C_s(y)
\end{equation}
where
\begin{equation}
C_s(y)\equiv \sumls b_\ell \cos(c_\ell y), \quad
S_s(y)\equiv \sumls b_\ell \sin(c_\ell y).
\end{equation}
The stability condition therefore becomes
\[
|R(\mathbf{i}y)|^2 = 1-2yS_s(y) + y^2 \left( C_s^2(y)+S_s^2(y) \right) \leq 1.
\]
Such condition can be written in the form
\begin{equation}
    y F_s(y) \geq 0, \quad {\rm where} \quad F_s(y)\equiv S_s(y) - \frac{1}{2}y\left(C_s^2(y)+S_s^2(y)\right).
   \label{eq:stab_inequal}
\end{equation}
The problem of finding quadrature formulas with the widest stability region
can be stated as: determine the coefficients
$\mathbf{b}=(b_1, \ldots, b_s)$ and $\mathbf{c}=(c_1, \ldots, c_s)$ so that the interval in which
\rf{eq:stab_inequal} is satisfied is the widest.

Rather than directly solving this optimization problem, we consider a particular case of
quadrature formulas, i.e. those for which the nodes are symmetrically
located with respect to point $1/2$ in the interval $[0,1]$.
Among such formulas we restrict to the case in which $s$ even, as the schemes are observed to be unstable for odd $s$, see Fig.~\ref{fig:glstab}.

Let us denote by $\tc_\ell = 1-2 c_\ell, \ell = 1, \ldots, s$.
Then $c_\ell  = (1-\tc_\ell)/2$.
Since the nodes are symmetric and the quadrature formula is
interpolatory, we have
\begin{equation}
   \tc_\ell = -\tc_{s-\ell+1}, \quad b_\ell = b_{s-\ell+1}.
   \label{eq:symmetry}
\end{equation}
The absolute stability function $R(\mathbf{i}y)$ can then be written,
after simple manipulations
\[
   R(\mathbf{i}y) = 1-2y \sin(y/2) \tC_s(y) + 2 \mathbf{i}y \cos(y/2)\tC_s(y),
\]
where
\[
   \tC_s(y)\equiv \suml{s/2}b_\ell \cos(\tc_\ell y/2),
\]
leading to
\begin{equation}
   |R(\mathbf{i}y)|^2 = 1-4y \sin(y/2)\tC_s(y) + 4 y^2 \tC_s^2(y).
   \label{eq:R2}
\end{equation}
The function $F_s(y)$ can then be written,
after simple manipulations
\begin{equation}
   F_s(y) = 2 \tC_s(y) \left(\sin(y/2) - y \tC_s(y)\right).
   \label{eq:F_s}
\end{equation}
Then the stability condition \rf{eq:stab_inequal} becomes
\[
 \tC_s(y) \left(\sin(y/2) - y \tC_s(y)\right) \geq 0.
\]
Because function $F_s$ contains the product between two factors, the condition to ensure that the function does not change sign at roots is that the two factors vanish simultaneously at simple roots, therefore $\tC_s(y)$ has to vanish also at the same
points $y_k > 0$ at which $\sin(y/2) - y \tC_s(y)=0$. There is no need to impose that
$\tC_s$ vanishes at the origin, since, because of symmetry, $yF_s(y)$ does not change sign at the origin.

In order to determine the coefficients that define the quadrature formula for maximizing the stability interval on imaginary axis, we proceed as follows.
Because of the symmetry constraints \rf{eq:symmetry}, we have to find
$s$ coefficients, i.e. $b_1,\ldots,b_{s/2}$ and $\tc_1,\ldots,\tc_{s/2}$ by imposing a total of $s$
conditions. Such conditions will be a balance between accuracy and
stability.
If we want that the quadrature formulas have degree of precision
$s-1$,
i.e. if we want that they are exact on polynomials of degree less or
equal to $s-1$, we have to impose
\begin{equation}
   \frac{1}{2} \int_{-1}^1 \zeta^{2k} \, dz =  \sumls b_\ell (\tc_\ell)^{2k},
   \quad k=0,\ldots,s/2-1.
   \label{eq:order1}
\end{equation}
We only impose the condition for even polynomials, since odd polynomials are automatically satisfied because of
symmetry.
The condition that $\tC_s(y)$ vanishes when $\sin(y/2)-y\tC_s(y)$ vanishes becomes
\begin{equation}
   \tC_s(2\pi k) = 0, \quad k=1,\ldots,s/2.
   \label{eq:C=0}
\end{equation}
For $k=0$ the stability condition \rf{eq:C=0} is marginally satisfied since
$\tC_s(0) = \suml{s/2}b_\ell = 1/2$.
Eqs.~\rf{eq:order1} and \rf{eq:C=0} constitute a nonlinear set of
equations for the $s$ coefficients $b_1,\ldots,b_{s/2}$ and
$\tc_1,\ldots,\tc_{s/2}$. Because the equations are nonlinear, we have
to resort to Newton's method for its solution. In practice, for
large values of $s$, it is hard to find an initial guess which lies in
the convergence basin of Newton's method. We had to resort to
a relaxed version of Newton's method, coupled with continuation
techniques, in order to solve the system.

\begin{figure}
\centering
\includegraphics[width=3.2in]{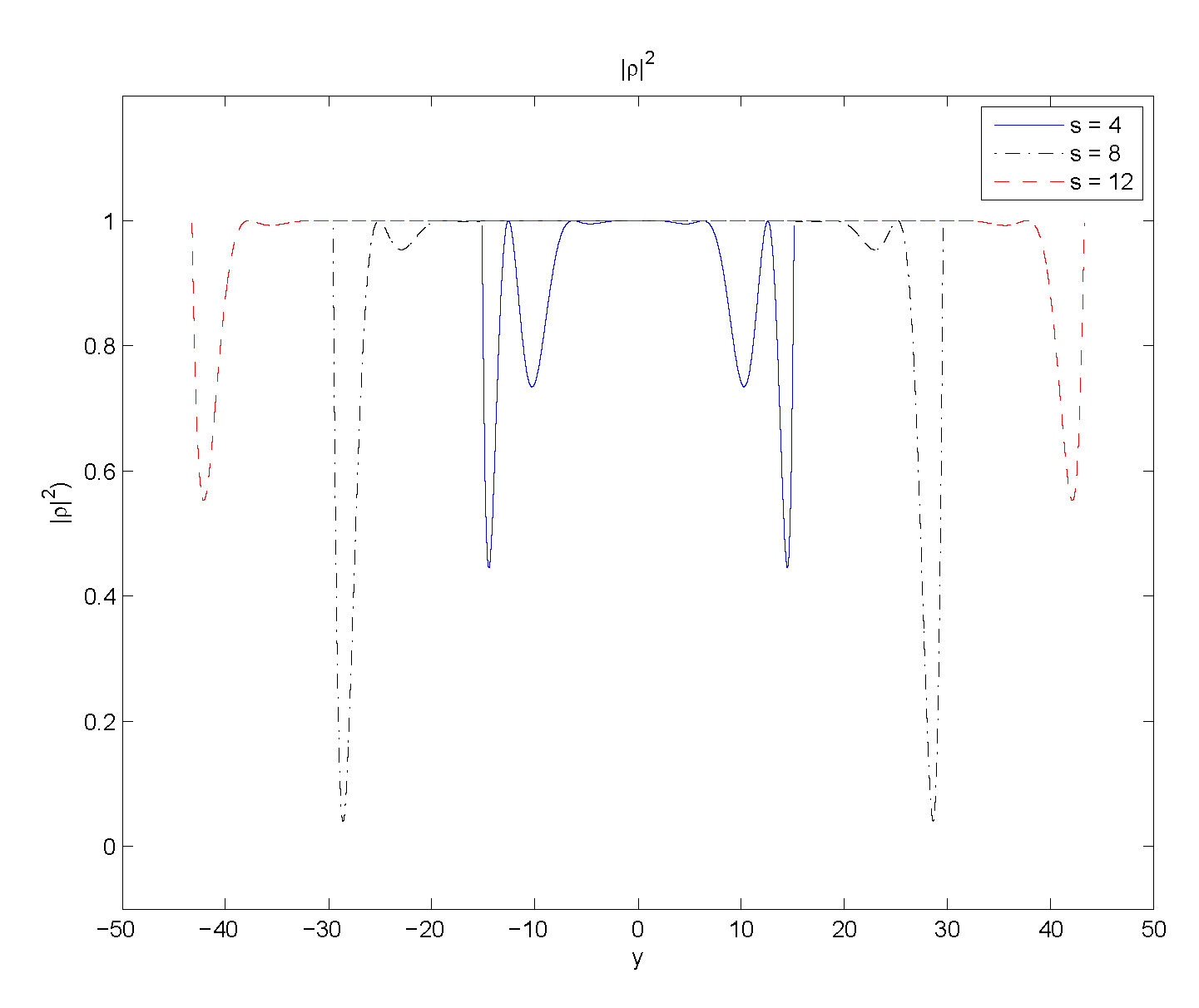}
\includegraphics[width=3.0in]{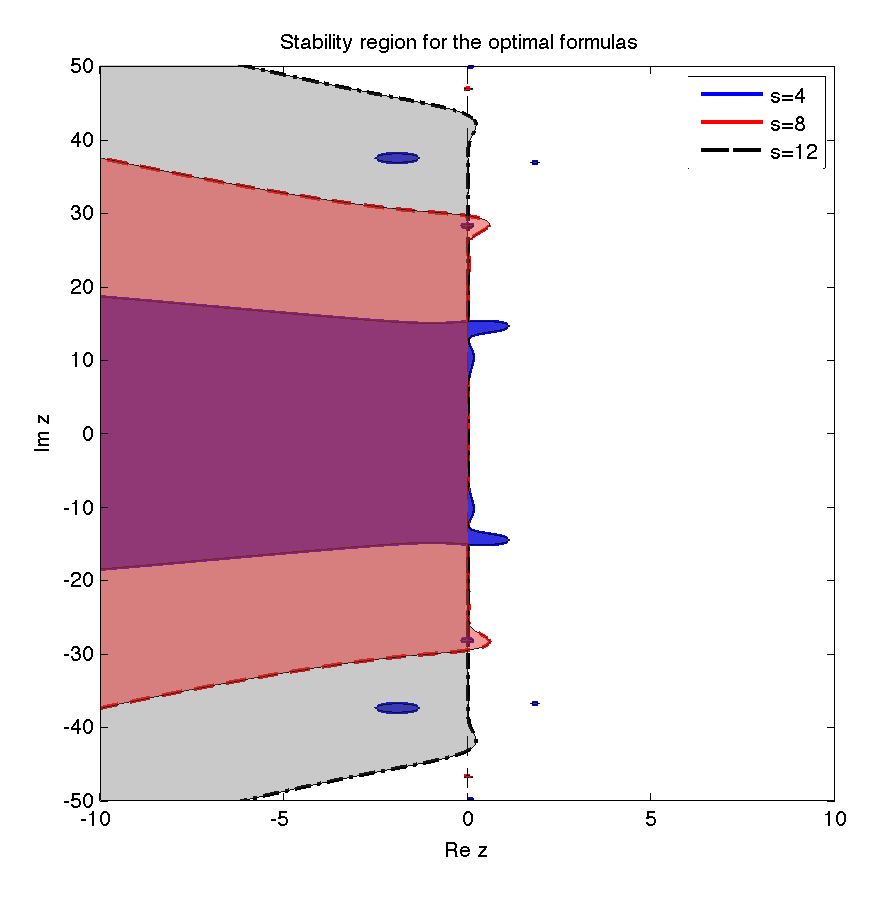}
\caption{Left: plot of $|R(\mathbf{i}y)|^2-1$ for the symmetric quadrature formulas with $s = 4, 8, 12$. Right: plot of stability regions for the corresponding quadrature formulas. These formulas show a wide stability interval on the imaginary axis, thus allowing, in principle, large CFL numbers.}
\label{fig:stab:check1}
\end{figure}

\begin{figure}
\centering
\includegraphics[width = 3.2in]{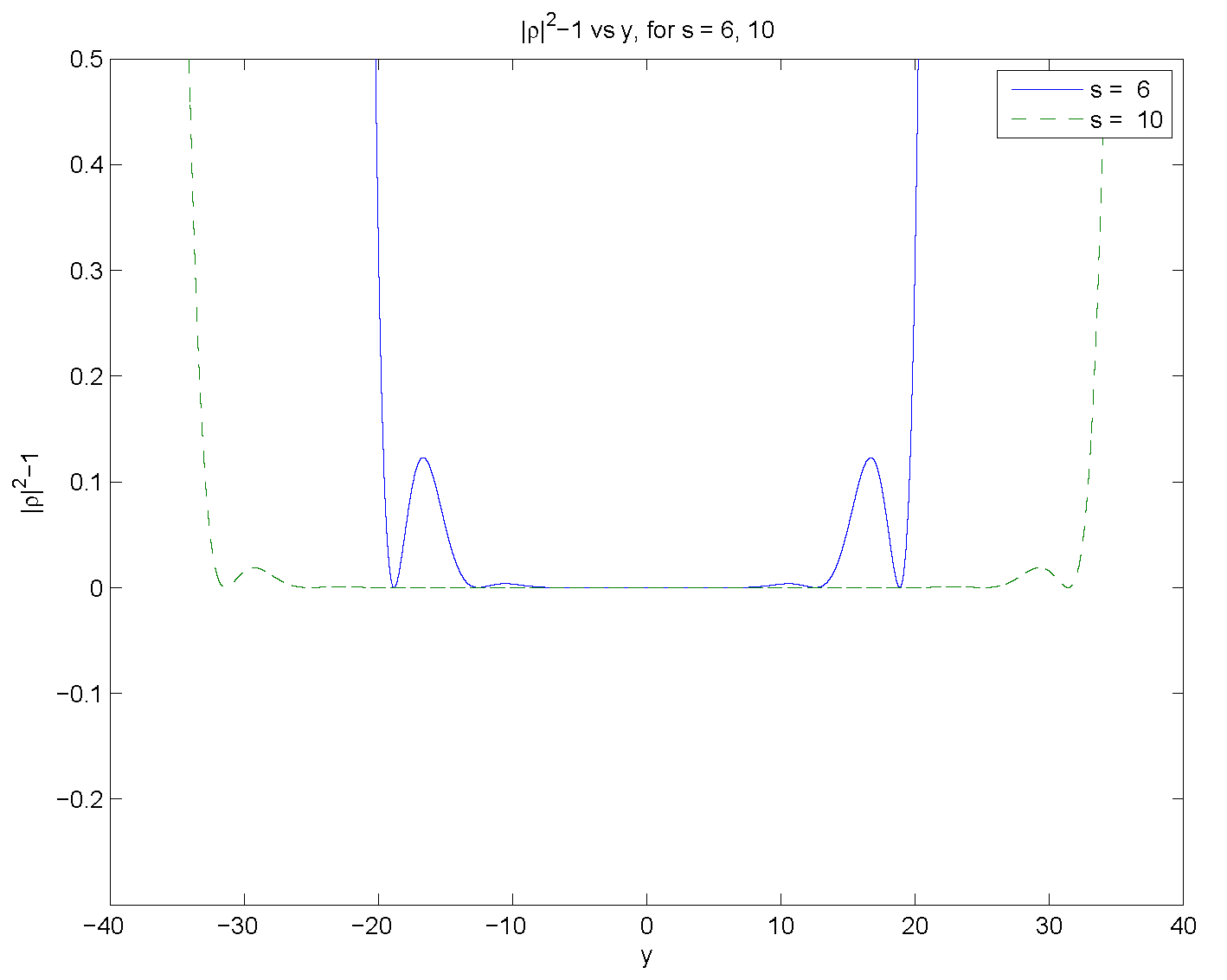}
\caption{Plot of $|R(\mathbf{i}y)|^2-1$ for the symmetric quadrature formulas
  with $s = 6, 10$. Such formulas are not stable because the
  stability region do not contain a portion of the imaginary axis.}
\label{fig:stab:check2}
\end{figure}

We numerically compute nodes and weights for $s=2, 4, 6, 8, 10, 12$ and check {\em a posteriori} whether the stability condition is actually
satisfied. The following phenomena are observed:
\bit
\item $s=2$: the quadrature nodes and weights are consistent with those in the two-point Gauss-Legendre formula.
\item $s = 4, 8, 12$: In Fig. \ref{fig:stab:check1} we plot the functions $|R_s(\mathbf{i}y)|^2-1$ (left panel) and the corresponding stability regions in the complex plane (right panel) for $s = 4, 8, 12$. A wide interval with stability on the imaginary axis is shown. We report the coefficients in Table \ref{tab:coeff}. Only $s/2$ coefficients are reported, since the other satisfy the symmetry relation \rf{eq:symmetry}. We also report the maximum CFL number $a^* =y^*/\pi$ (see eq. \eqref{eq: cfl_ode}).
\item $s=6$ and $s=10$: In Fig. \ref{fig:stab:check2}, we plot the functions $|R_s(\mathbf{i}y)|^2-1$ for $s = 6$ and $s = 10$. It is observed that $R_s(\mathbf{i}y)|^2\ge1$ for any interval containing the origin, i.e. these two quadrature formulas are not stable.
\eit

\begin{table}
\begin{center}
\caption{Weights and nodes of accurate and stable quadrature
  formulas. Each formula is exact for polynomials of degree not
  greater than $s-1$. The maximum CFL number $a^*$ that guarantees
  stability in the theoretical Fourier analysis is reported.}
  \bigskip
\begin{tabular}{|c|c|c|} \hline
\multicolumn{3}{|c|}{$s=4,  a^*= 4.8125674352016$} \\ \hline
1 & 0.199889211759008 & 0.083205952308564 \\
2 & 0.300110788240992 & 0.347904700949451 \\ \hline
\multicolumn{3}{|c|}{$s=8,  a^*= 9.4130380474585$} \\ \hline
1 & 0.058702317190867 & 0.023248965963790 \\
2 & 0.119923212650690 & 0.114686793929813 \\
3 & 0.154113350301760 & 0.253867587586135 \\
4 & 0.167261119856682 & 0.415892817555109 \\ \hline
\multicolumn{3}{|c|}{$s=12, a^*=13.7671988660496$} \\ \hline
1 & 0.027182888487959 & 0.010668025829619 \\
2 & 0.059633412276882 & 0.054560771376909 \\
3 & 0.084799522112170 & 0.127471263371368 \\
4 & 0.101625491473440 & 0.221353922812027 \\
5 & 0.111259037829236 & 0.328318059665840 \\
6 & 0.115499647820313 & 0.442082833046309 \\ \hline
\end{tabular}
\label{tab:coeff}
\end{center}
\end{table}

The stability regions for the quadrature formulas obtained for
$s=4,8,12$ reported in Table \ref{tab:coeff} are computed under the
assumption that one considers the exact space dependence of the
Fourier mode, so that the only error is in time integration. In
reality there are several other causes of errors, that may affect the
stability region of the quadrature. In the next section,
we take spatial discretization into account and quantify the corresponding stability interval.

%
%

%% file: pde_stab.tex
\section{Spatial discretization.}
\label{ssec: spat}
\setcounter{equation}{0}
\setcounter{figure}{0}
\setcounter{table}{0}

There are two spatial discretization processes in the scheme. One is the WENO interpolation in approximating $f(x_j, t^n+c_\ell \Delta t) = f(x_j-c_\ell \Delta t/\Delta x, t^n)$ from neighboring grid point values $\{f^n_j\}_{j=1}^{n_x}$. The other is the WENO reconstruction in obtaining numerical fluxes $\hat{F}_{j+\frac12}$ in \eqref{eq: cons_update} from $\{\mathcal{F}_j\}_{j=1}^{nx}$. In this paper, we consider the following two classes of spatial discretizations.
\begin{itemize}
\item Odd order approximations. For the linear equation \eqref{eq:scalar2}, we use a right-biased stencil to approximate $f(x_j-c_\ell \Delta t/\Delta x, t^n)$ and use a left-biased stencil for reconstructing the flux $\hat{F}_{j+\frac12}$. For example, for a first order scheme with $ \Delta t/\Delta x<1$, $f(x_j-c_\ell \Delta t/\Delta x, t^n)$ is approximated from the {\em interpolation} stencil $\{f_j\}$
and the numerical flux $\hat{F}_{j+\frac12}$ is approximated from the {\em reconstruction} stencil $\{\mathcal{F}_j\}$. With such stencil arrangement, the SL scheme is reduced to a first order upwind scheme when $ \Delta t/\Delta x<1$,
\[
f^{n+1}_j = f^n_j - \Delta t/\Delta x (f^n_j - f^n_{j-1}).
\]
Third, fifth, seventh and ninth order schemes can be constructed by including one, two, three, four more points symmetrically from left and from right, respectively, in the interpolation and reconstruction stencils. We list them as follows.
\beq
\begin{array}{lll}
&\mbox{Third order}: & \{f_{j-1}, f_j, f_{j+1}\}, \quad \{\mathcal{F}_{j-1}, \mathcal{F}_j, \mathcal{F}_{j+1}\}.\\[2mm]
&\mbox{Fifth order}: & \{f_{j-2}, f_{j-1}, f_j, f_{j+1}, f_{j+2}\}, \quad \{\mathcal{F}_{j-2}, \mathcal{F}_{j-1}, \mathcal{F}_j, \mathcal{F}_{j+1}, \mathcal{F}_{j+2}\}.\\[2mm]
&\mbox{Seventh order}: & \{f_{j-3}, f_{j-2}, f_{j-1}, f_j, f_{j+1}, f_{j+2}, f_{j+3}\}, \\
& & \{\mathcal{F}_{j-3}, \mathcal{F}_{j-2}, \mathcal{F}_{j-1}, \mathcal{F}_j, \mathcal{F}_{j+1}, \mathcal{F}_{j+2}, \mathcal{F}_{j+3}\}.\\[2mm]
&\mbox{Ninth order}: & \{f_{j-4}, f_{j-3}, f_{j-2}, f_{j-1}, f_j, f_{j+1}, f_{j+2}, f_{j+3}, f_{j+4}\}, \\
& & \{\mathcal{F}_{j-4}, \mathcal{F}_{j-3}, \mathcal{F}_{j-2}, \mathcal{F}_{j-1}, \mathcal{F}_j, \mathcal{F}_{j+1}, \mathcal{F}_{j+2}, \mathcal{F}_{j+3}, \mathcal{F}_{j+4}\}.
\end{array}
\notag
\eeq
\item Even order approximations. For the linear equation \eqref{eq:scalar2}, we use symmetric stencils to approximate $f(x_j-c_\ell \Delta t/\Delta x, t^n)$ by interpolation and to approximate $\hat{F}_{j+\frac12}$ by reconstruction. For example, for a second order scheme with $ \Delta t/\Delta x<1$, $f(x_j-c_\ell \Delta t/\Delta x, t^n)$ is approximated from the {\em interpolation} stencil $\{f_{j-1}, f_j\}$
and the numerical flux $\hat{F}_{j+\frac12}$ is approximated from the {\em reconstruction} stencil $\{\mathcal{F}_j, \mathcal{F}_{j+1}\}$. Fourth, sixth and eighth order schemes can be constructed by including one, two, three more points symmetrically from left and from right, respectively, in the interpolation and reconstruction stencils. We list them as follows.
\beq
\begin{array}{lll}
&\mbox{Fourth order}: & \{f_{j-2}, f_{j-1}, f_j, f_{j+1}\}, \quad \{\mathcal{F}_{j-1}, \mathcal{F}_j, \mathcal{F}_{j+1}, \mathcal{F}_{j+2}\}.\\[2mm]
&\mbox{Sixth order}: & \{f_{j-3}, f_{j-2}, f_{j-1}, f_j, f_{j+1}, f_{j+2}\}, \\
& & \{\mathcal{F}_{j-2}, \mathcal{F}_{j-1}, \mathcal{F}_j, \mathcal{F}_{j+1}, \mathcal{F}_{j+2}, \mathcal{F}_{j+3}\}.\\[2mm]
&\mbox{Eighth order}: & \{f_{j-4}, f_{j-3}, f_{j-2}, f_{j-1}, f_j, f_{j+1}, f_{j+2}, f_{j+3}\}, \\
& & \{\mathcal{F}_{j-3}, \mathcal{F}_{j-2}, \mathcal{F}_{j-1}, \mathcal{F}_j, \mathcal{F}_{j+1}, \mathcal{F}_{j+2}, \mathcal{F}_{j+3}, \mathcal{F}_{j+4}\}.
\end{array}
\notag
\eeq
\end{itemize}
\begin{rem}
We follow the same principle in the interpolation and reconstruction procedures in more general settings, for example the situation when the time stepping size is greater than the CFL restriction, i.e $\Delta t/\Delta x \ge 1$ for eq. \eqref{eq:scalar2}. For general high dimensional problems, e.g. the Vlasov equation, similar procedures can be applied in a truly multi-dimensional fashion.
\end{rem}

To access the stability property of the conservative method, we perform Fourier analysis via the linear equation \eqref{eq:scalar2} with $x\in [0, 2\pi]$ and periodic boundary condition. In particular, we make the ansatz $f^n_j = \hat{f}^n e^{\mathbf{i} j \xi}$ with $\mathbf{i} = \sqrt{-1}$ and $\xi \in [0, 2\pi]$. Plugging the ansatz into the SL conservative scheme as described in Section~\ref{sec: cons}, we obtain $\hat{f}^{n+1}(\xi) = Q_{\lambda}(\xi) \hat{f}^n (\xi)$ with $Q_{\lambda}(\xi)$ being the amplification factor for the Fourier mode associated with $\xi$ and $\lambda = \frac{\Delta t}{\Delta x}$. To ensure linear stability, it is sufficient to have
\beq
\label{eq: linear_stab}
|Q_\lambda (\xi)|\le 1, \quad \forall \xi \in [0, 2\pi], \quad \forall\lambda \in [0, \lambda^\star], \quad \mbox{for some $\lambda^\star$}.
\eeq
We seek for $\lambda^\star$ by numerically checking the inequality \eqref{eq: linear_stab} for $100$ discretized grid points on $\xi \in [0, 2\pi]$,
and by gradually increasing $\lambda$ with a step size of $0.01$ starting from $\lambda =0$. Taking the machine precision into account in our implementation, we check the inequality $|Q_\lambda (\xi)|\le 1 + 10^{-11}$ instead. We tabulate such $\lambda^\star$ in Table~\ref{tab41} for different quadrature formulas as discussed in Section~\ref{ssec: temp} and with different choices of spatial interpolation and reconstruction stencils with odd and even order respectively. One can observe that the second order trapezoidal rule and the fourth order GL2 perform much better than the mid-point rule in terms of stability, especially when the orders for spatial approximations are high. The time stepping sizes allowed for stability of fully discretized schemes with $s=4, 8, 12$ are observed to be much less than the one provided by ODE stability analysis in the previous section.

\begin{figure}
\begin{center}
\includegraphics[height=2.5in,width=3.0in]{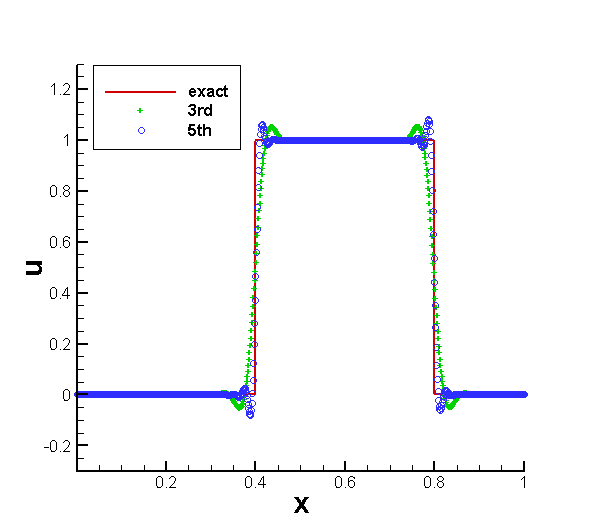}
\includegraphics[height=2.5in,width=3.0in]{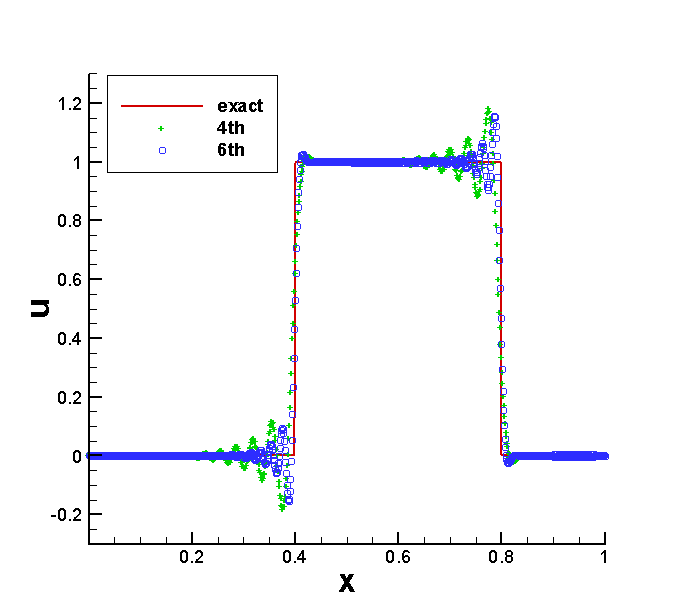}
\end{center}
\caption{Numerical solution for the linear advection equation $u_t + u_x = 0$ with an initial step function
at $T=100.1$. Left: 3rd order with $CFL=1.22$ and 5th order with $CFL=1.19$; Right: 4th order and 6th order with $CFL=1.84$. $N=800$. Here linear interpolation and linear reconstruction are used without WENO.}
\label{fig:test}
\end{figure}

In the following, we take the linear advection equation $u_t + u_x = 0$ with a smooth initial function $\sin(2\pi x)$ on the domain $[0,1]$, to test the CFL bounds in Table \ref{tab41}. Here for better illustration, only linear interpolation and linear reconstruction are used. We consider schemes that couple GL2 for temporal integration with third and fourth order spatial approximations. Errors and orders of convergence at a final integration time $T=100.1$ are recorded in Table \ref{tab42}. Clear third order and fourth order spatial accuracy are observed at the corresponding upper bounds for CFL ($1.22$ for third order and $1.84$ for fourth order as in Table \ref{tab41}.) The code will blow up with the CFL increased by $0.01$ at the corresponding time, which confirms the validity of the CFL bounds in the table. We have similar observations for other orders of schemes, but omit to present them to save space. Although even order schemes comparatively have larger CFL bounds than odd order ones, for solutions with discontinuities, we can observe that odd order schemes with upwind mechanism can resolve the discontinuities better. We present numerical solutions of our schemes with linear weights for advecting a step function in Fig. \ref{fig:test}. 
Due to the above considerations, we use the scheme with the 5th order spatial approximation and with two-point Gaussian rule for temporal integration in the following numerical sections.

\begin{table}
\begin{center}
\caption{Upper bounds of CFL for FD SL scheme with odd and even order interpolation and reconstruction.
The amplification factor is bounded by $1+10^{-11}$. $N=100$.}
\bigskip
\begin{tabular}{|c|c|c|c|c|c|c|} \hline
temp/spatial& 1st& 3rd  & 5th & 7th & 9th & exact \\\hline
mid-point &    1.00 &    1.00 &    0.14 &    0.04 &    0.02 & 0.00 \\\hline
trapezoid &    1.99 &    1.68 &    1.52 &    1.44 &    1.38 & 1.00\\\hline
Simpson   &    1.33 &    1.50 &    1.35 &    0.71 &    0.37 & 0.00 \\\hline
GL2       &    1.00 &    1.22 &    1.19 &    1.16 &    1.15 & 1.72  \\\hline
s=4       &    1.00 &    1.37 &    1.27 &    1.22 &    1.19 & 4.81 \\\hline
s=8       &    1.00 &    1.35 &    1.26 &    1.21 &    1.18 & 9.41\\\hline
s=12      &    1.00 &    1.37 &    1.25 &    1.21 &    1.18 & 13.76\\\hline \hline
temp/spatial& 2nd & 4th  & 6th & 8th & 10th & exact \\\hline
mid-point &    2.00 &    0.04 &    0.01 &    0.00 &    0.00 &0.00 \\\hline
trapezoid &    1.29 &    1.26 &    1.24 &    1.22 &    1.20 &1.00 \\\hline
Simpson   &    3.00 &    2.91 &    0.83 &    0.34 &    0.20 &0.00 \\\hline
GL2       &    1.85 &    1.84 &    1.84 &    1.83 &    1.83 &1.72 \\\hline
s=4       &    1.96 &    1.97 &    1.98 &    1.98 &    1.98 &4.81 \\\hline
s=8       &    1.99 &    1.99 &    1.99 &    1.99 &    1.99 &9.41 \\\hline
s=12      &    1.99 &    1.99 &    1.99 &    1.99 &    2.00 &13.76 \\\hline
\end{tabular}
\label{tab41}
\end{center}
\end{table}

\begin{table}
\begin{center}
\caption{Accuracy test of the linear advection equation $u_t + u_x = 0$ with the initial function $\sin(2\pi x)$ for the 3rd order scheme with $CFL=1.22$ at $T=100.1$ and 4th order with $CFL=1.84$ at $T=1001.1$.}
\bigskip
\begin{tabular}{|c|c|c|c|c|c|} \hline
Scheme &   N  &  $L^1$ error & order   & $L^\infty$ error & order \\  \hline
\multirow{4}{*}{3rd order} &  240 &     4.12E-04 &       --&     6.47E-04 &       --  \\  \cline{2-6}
 &  480 &     5.15E-05 &     3.00&     8.09E-05 &     3.00  \\  \cline{2-6}
 &  960 &     6.44E-06 &     3.00&     1.01E-05 &     3.00  \\  \cline{2-6}
 & 1920 &     8.04E-07 &     3.00&     1.26E-06 &     3.00  \\  \hline
\multirow{4}{*}{4th order}&  120 &     1.76E-03 &       --&     2.77E-03 &       --  \\  \cline{2-6}
 &  240 &     1.10E-04 &     4.00&     1.73E-04 &     4.00  \\  \cline{2-6}
 &  480 &     6.89E-06 &     4.00&     1.08E-05 &     4.00  \\  \cline{2-6}
 &  960 &     4.76E-07 &     3.85&     7.48E-07 &     3.85  \\  \hline
\end{tabular}
\label{tab42}
\end{center}
\end{table}

%% file: numerics.tex
\section{Numerical tests on 2D linear passive-transport problems}
\label{sec5}
\setcounter{equation}{0}
\setcounter{figure}{0}
\setcounter{table}{0}

In this section, the conservative truly multi-dimensional SL scheme will be tested for passive transport equations, such as linear advection, rotation and swirling deformation. Since the velocity of the field is given a priori, characteristics can be traced by a high order Runge-Kutta ODE integrator.

In this and next sections, we use $5$th order spatial approximations with WENO (i.e. WENO interpolation and WENO reconstruction) for evaluating flux functions in \eqref{eq: cons_update}.
We use GL2 for temporal integration, while characteristics are traced back in time by Runge-Kutta to locate feet of characteristics. For a general two dimensional problem $u_t+f(u)_x + g(u)_y = 0$, the time step is taken as
\[
\Delta t = CFL/(a/\Delta x+b/\Delta y),
\]
where $a=\max|f'(u)|$ and $b=\max|g'(u)|$. From Table \ref{tab41}, the CFL number is $1.22$ for a 3rd order spatial discretization and $1.19$ for the $5$th order. In the following, we take $CFL=1.15$ without specification.



\begin{exa}
We first test our problem for the linear equation $u_t+u_x+u_y=0$ with initial condition
$u(x,y,0)=\sin(x)\sin(y)$. The exact solution is $u(x,y,t)=\sin(x-t)\sin(y-t)$.
For this example, the roots of characteristics are located exactly.
Table \ref{tab1} and Table \ref{tab1t} presents spatial and temporal order of convergence of the proposed scheme.
Both $5$-th order spatial accuracy and $4$-th order temporal accuracy from GL2 can be observed.

\begin{table}[h]
\begin{center}
\label{tab1}
\caption{{Errors and orders for the linear equation in space. $T=1.2$. $CFL=1.15$.}}
\bigskip
\begin{tabular}{|c | c|c|c| c|c|}
\hline
\cline{1-5} $N_x \times N_y$  &{$20\times 20$} &{$40\times 40$} &{$60\times 60$}&{$80\times 80$}&{$100\times100$}  \\
\hline
\cline{1-5}  $L^1$ error &2.76E-4&8.38E-6 &1.11E-6&2.64E-7&8.68E-8\\
\hline
 order &--&5.04&4.99&4.98&4.99\\
\hline
\end{tabular}
\bigskip
\caption{{Errors and orders for the linear equation in time. $N_x=N_y=200$. $T=1$.}}
\bigskip
\begin{tabular}{|c | c|c|c| c|c|}
\hline
\cline{1-5} $CFL$  &{$1.1$} &{$1.0$} &{$0.9$}&{$0.8$}&{$0.7$}  \\
\hline
\cline{1-5}  $L^1$ error &3.34E-9&2.27E-9 &1.49E-9&9.29E-10&5.46E-10\\
\hline
 order &--&4.07&3.97&4.02&3.98\\
\hline
\end{tabular}
\label{tab1t}
\end{center}
\end{table}


\end{exa}

\begin{exa}
Now we consider two problems defined on the domain $[-\pi,\pi]^2$.
One is the rigid body rotating problem
\[
u_t - y u_x + x u_y =0,
\]
the other is the swirling deformation flow problem
\[
u_t - \left(\cos^2\left(\frac{x}{2}\right)\sin(y)g(t)u\right)_x + \left(\sin(x)\cos^2\left(\frac{y}{2}\right)g(t)u\right)_y = 0,
\]
with $g(t)=\cos(\pi t/T)\pi$. Both have the initial condition, which includes a slotted disk, a cone as well as a smooth hump, see Fig. \ref{fig1} (top) and Fig. \ref{fig2} (top).

For the rigid body rotating problem, its period is $2\pi$. In Fig. \ref{fig1}, we have shown
the results at a half period and one period. As we can see, the shape of the bodies are well preserved. For the swirling deformation flow problem, after a half period the bodies are deformed, but they regain its initial shape after one period, see Fig. \ref{fig2}.

\begin{figure}
\begin{center}
\includegraphics[height=2.5in,width=3.0in]{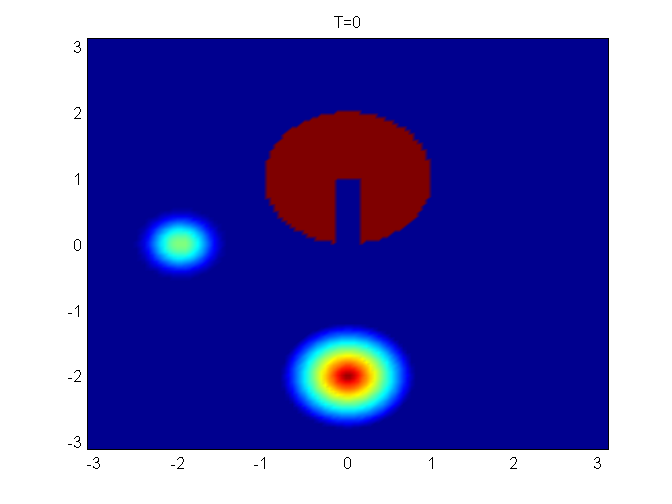}
\includegraphics[height=2.5in,width=3.0in]{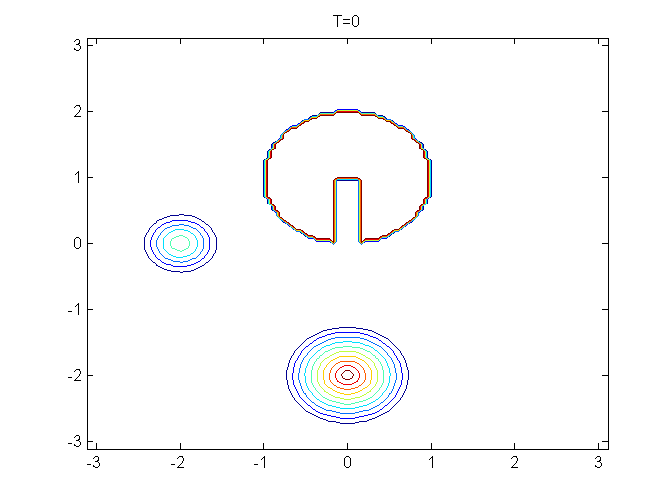}\\
\includegraphics[height=2.5in,width=3.0in]{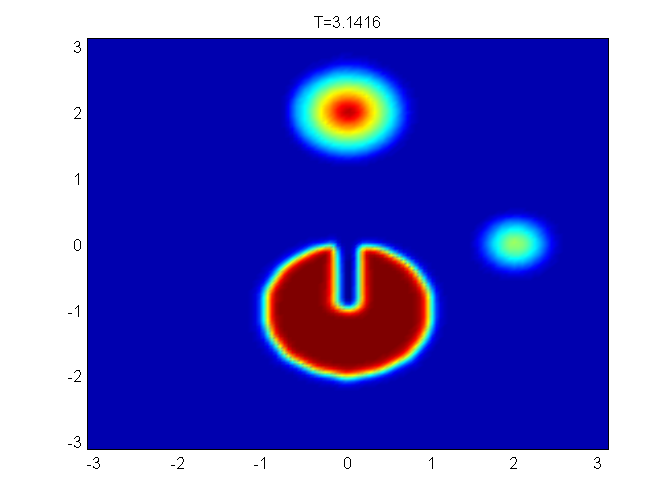}
\includegraphics[height=2.5in,width=3.0in]{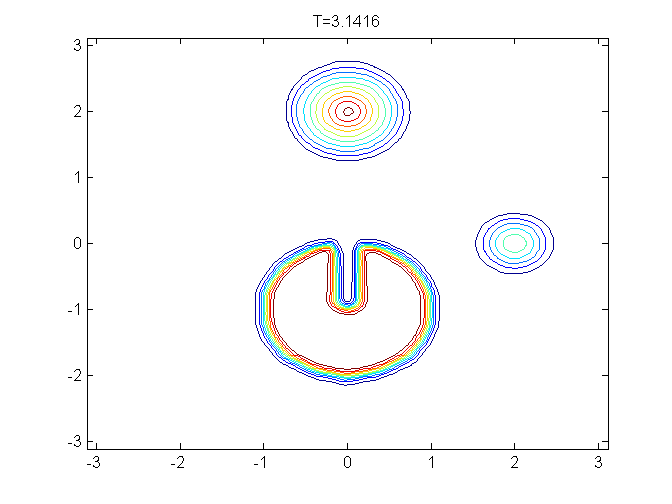}\\
\includegraphics[height=2.5in,width=3.0in]{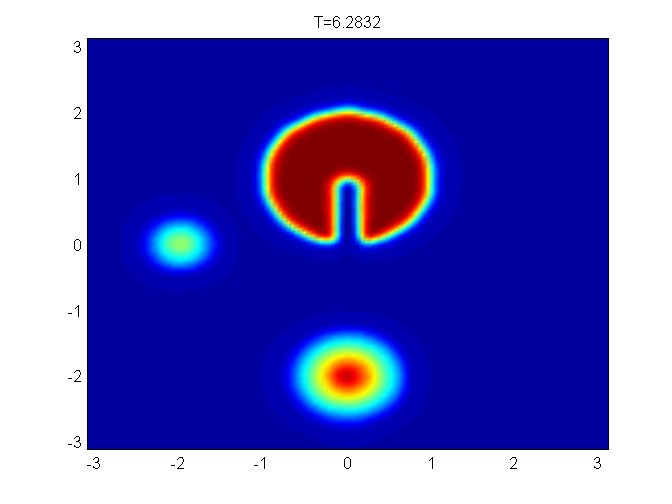}
\includegraphics[height=2.5in,width=3.0in]{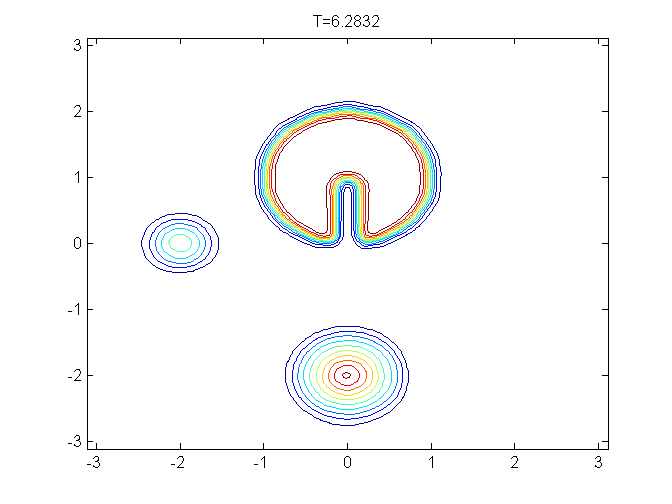}
\end{center}
\caption{Rigid body rotating problem. Mesh size: $128\times128$. Top: $T=0$; middle: $T=\pi$;
bottom: $T=2\pi$. Contour plots: 10 equally spaced lines.}
\label{fig1}
\end{figure}

\begin{figure}
\begin{center}
\includegraphics[height=2.5in,width=3.0in]{./pic/rigid/surface.png}
\includegraphics[height=2.5in,width=3.0in]{./pic/rigid/contour.png}\\
\includegraphics[height=2.5in,width=3.0in]{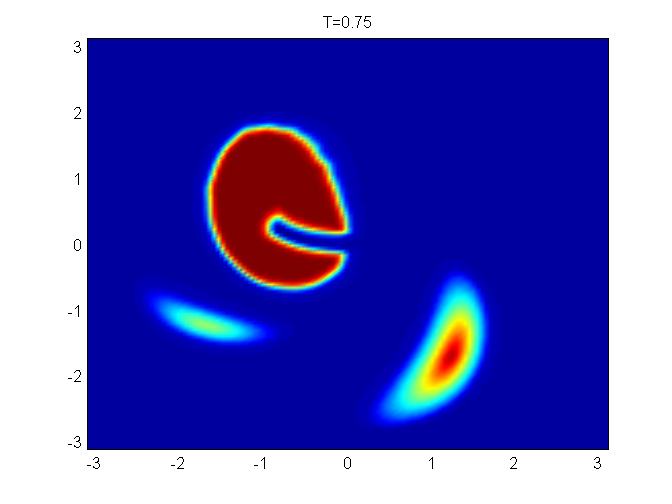}
\includegraphics[height=2.5in,width=3.0in]{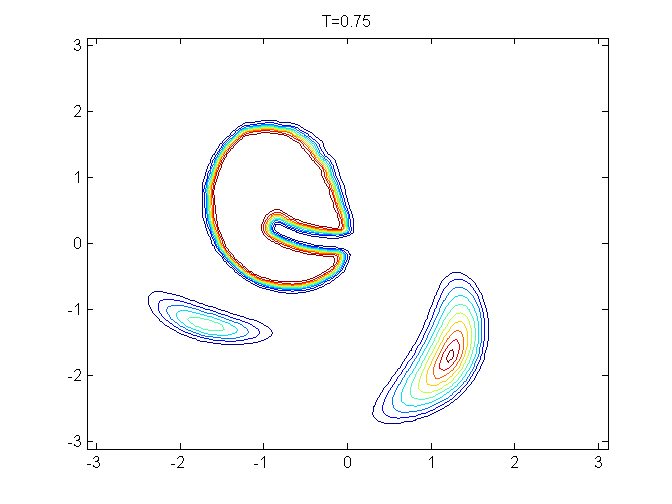}\\
\includegraphics[height=2.5in,width=3.0in]{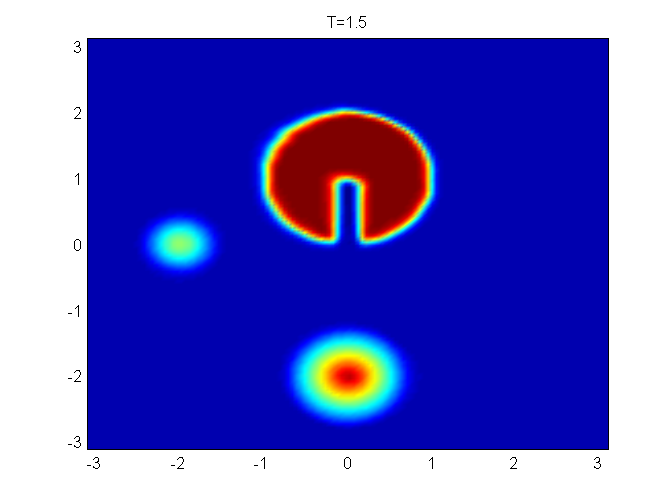}
\includegraphics[height=2.5in,width=3.0in]{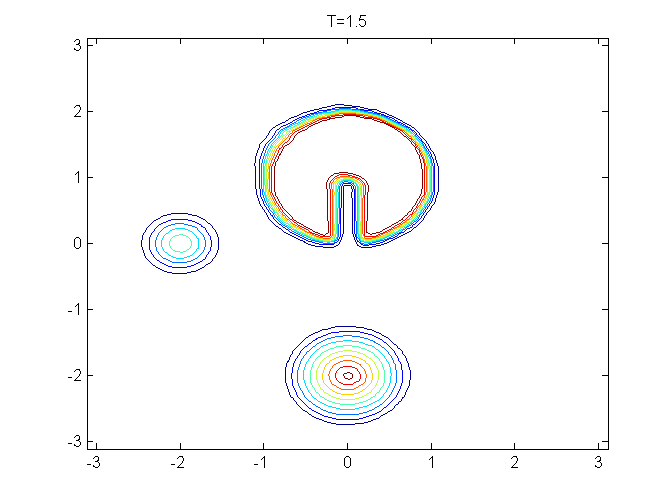}
\end{center}
\caption{Rigid body rotating problem. Mesh size: $128\times128$. Top: $T=0$; middle: $T=0.75$;
bottom: $T=1.5$. Contour plots: 10 equally spaced lines.}
\label{fig2}
\end{figure}
\end{exa}

\section{Numerical tests of nonlinear systems}
\label{sec6}
\setcounter{equation}{0}
\setcounter{figure}{0}
\setcounter{table}{0}

In this section, we test the conservative SL scheme on the nonlinear VP system, the guiding center Vlasov system and the incompressible Euler system in vorticity stream function formulation. Despite different application backgrounds, the latter two systems are indeed in almost the same mathematical formulation, only with different signs in the Poisson's equation.

\subsection{VP system}
Arising from collisionless plasma applications, the VP system
\begin{equation}
\frac{\partial f}{\partial t} + {\bf v} \cdot \nabla_{\bf x} f  +
\mathbf{E}({\bf x},t) \cdot \nabla_{\bf v} f = 0, \label{eq: vlasov}
\end{equation}
and
\begin{equation}
\mathbf{E}(\mathbf{x},t)=-\nabla_{\bf x}\phi(\mathbf{x},t),\quad
-\Delta_{\bf
x}\phi(\mathbf{x},t)=\rho(\mathbf{x},t)-1,\label{eq: poisson}
\end{equation}
describes the temporal evolution of the particle distribution function in six dimensional phase space. $f( {\bf x},{\bf v},t)$ is the
probability distribution function which describes the probability of
finding a particle with velocity $\bf{v}$ at position $\bf{x}$ at
time $t$, $\bf{E}$ is the electric field, and $\phi$ is the
self-consistent electrostatic potential. The probability
distribution function couples to the long range fields via the
charge density, $\rho(t,x) = \int_{\mathbb{R}^3} f(x,v,t)dv$,
where we take the limit of uniformly distributed infinitely massive
ions in the background. In this paper, we consider the VP system with
1-D in ${\bf x}$ and 1-D in ${\bf v}$. Periodic
boundary condition is imposed in x-direction, while zero boundary
condition is imposed in v-direction. The equations for tracking characteristics are
\beq
\frac{dx}{dt} = v, \quad \frac{dv}{dt} = E,
\eeq
where $E$ nonlinearly depends on $f$ via the Poisson system \eqref{eq: poisson}. To locate the foot of characteristics accurately, we apply the high order procedure proposed in \cite{qiu_russo_2016}.

Next we recall several norms in the VP system below, which should remain constant in time.
\begin{enumerate}
\item Mass:
\[
\text{Mass}=\int_v\int_xf(x,v,t)dxdv.
\]
\item $L^p$ norm $1\leq p<\infty$:
\begin{equation}
\|f\|_p=\left(\int_v\int_x|f(x,v,t)|^pdxdv\right)^\frac1p.
\end{equation}
\item Energy:
\begin{equation}
\text{Energy}=\int_v\int_xf(x,v,t)v^2dxdv + \int_xE^2(x,t)dx,
\end{equation}
where $E(x,t)$ is the electric field.
\item Entropy:
\begin{equation}
\text{Entropy}=\int_v\int_xf(x,v,t)\log(f(x,v,t))dxdv.
\end{equation}
\end{enumerate}
Tracking relative deviations of these quantities numerically will be
a good measure of the quality of numerical schemes. The relative
deviation is defined to be the deviation away from the corresponding
initial value divided by the magnitude of the initial value.
We also check the mass conservation over time $\int_v\int_x f(x,v,t)dxdv$, which is the same as the $L^1$ norm if $f$ is positive. However, since our scheme is not positivity preserving, the time evolution of the mass could be different from that of the $L^1$ norm due to the negative values appearing in numerical solutions.

In our numerical tests, we let the time step size $\Delta t = CFL \cdot \min(\Delta x/v_{max}, \Delta v/\max(E))$, where $CFL$ is specified as 1.15, and let $v_{max} = 6$ to minimize the error from truncating the domain in $v$-direction.

\begin{exa} (Weak Landau damping)
For the VP system, we first consider the weak Landau damping with the initial condition:
\beq
\label{landau}
f(t=0,x,v)=\f{1}{\sqrt{2\pi}}(1+\alpha \cos(k x))\exp(-\f{v^2}{2}),
\eeq
where $\alpha=0.01$ and $k=0.5$. The length of the domain in the x-direction is
$L=\f{2\pi}{k}$, which is similar in the following examples.
 In Fig. \ref{fig21}, we plot the time evolution of the electric field in $L^2$ norm and $L^\infty$ norm, the relative derivation of the discrete $L^1$ norm, $L^2$ norm, kinetic energy and entropy.

\begin{figure}
\centering
\includegraphics[width=3in,clip]{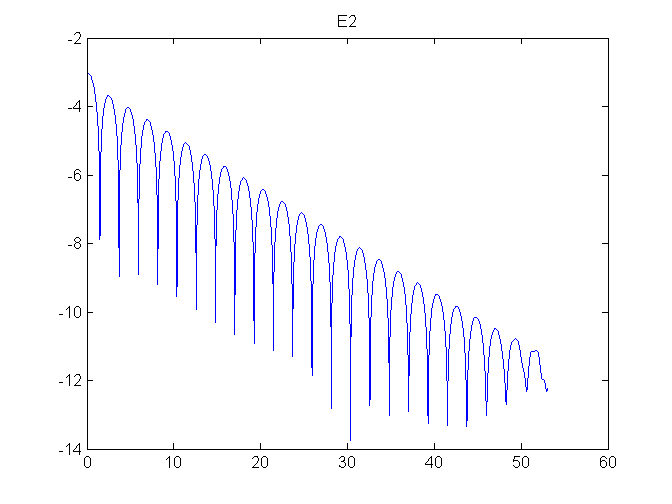},
\includegraphics[width=3in,clip]{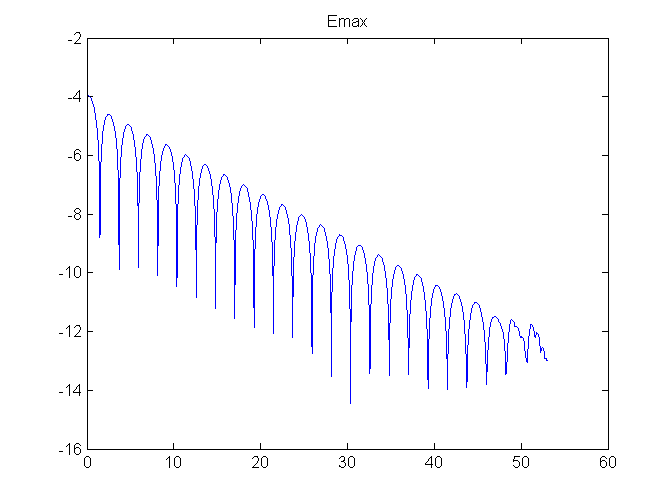} \\
\includegraphics[width=3in,clip]{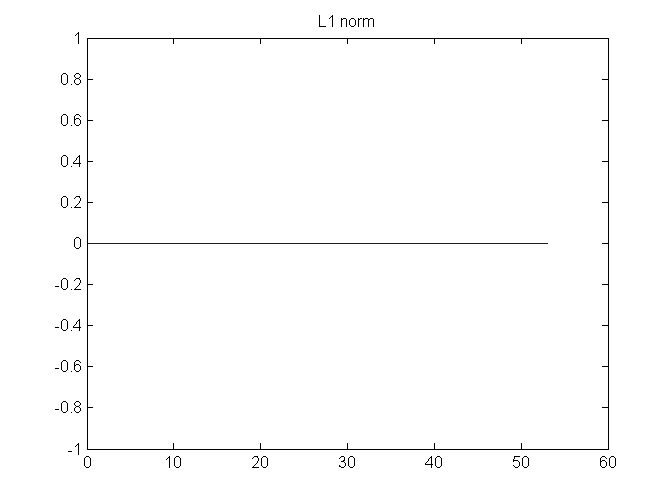},
\includegraphics[width=3in,clip]{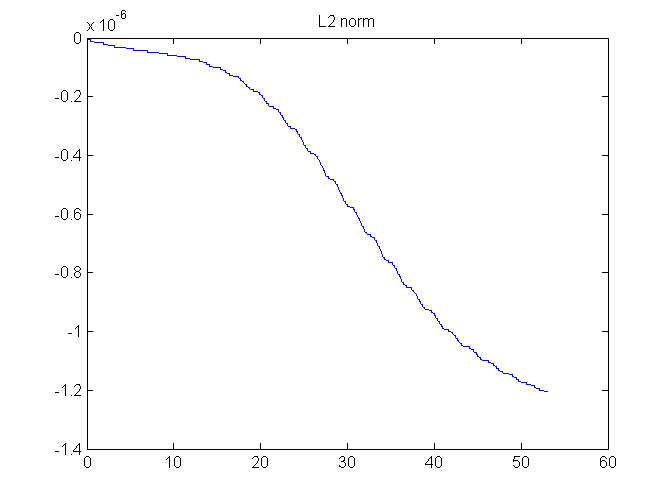}\\
\includegraphics[width=3in,clip]{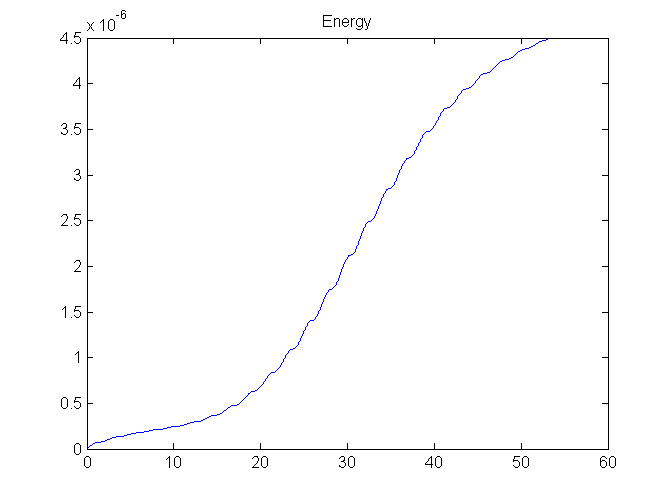},
\includegraphics[width=3in,clip]{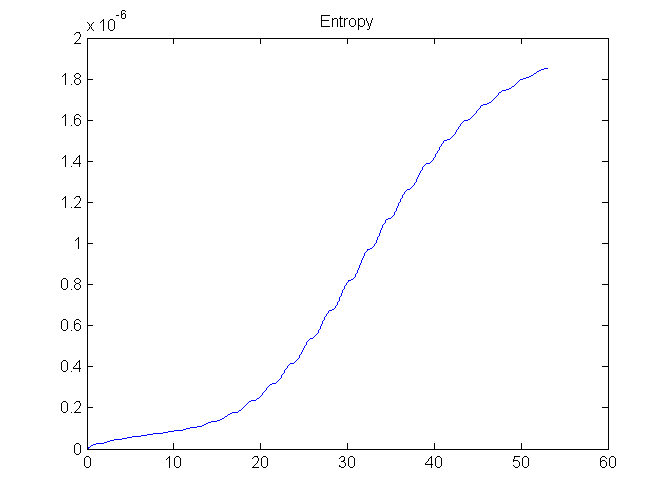}
\caption{Weak Landau damping. Time evolution of the electric field in $L^2$ norm and $L^\infty$ norm (top), discrete $L^1$ norm and $L^2$ norm (middle), kinetic
energy and entropy (bottom). Mesh: $128 \times 128$.}
\label{fig21}
\end{figure}
\end{exa}

\begin{exa} (Strong Landau damping)
The initial condition of strong Landau damping is still to be (\ref{landau}),
with $\alpha=0.5$ and $k=0.5$. Similarly in Fig. \ref{fig22}, we plot the time evolution
of electric field in $L^2$ norm and $L^\infty$ norm, the relative derivation of the discrete $L^1$ norm, $L^2$ norm, kinetic energy and entropy. The mass conservation is indicated by the bottom straight line in the $L^1$ norm figure.

\begin{figure}
\centering
\includegraphics[width=3in,clip]{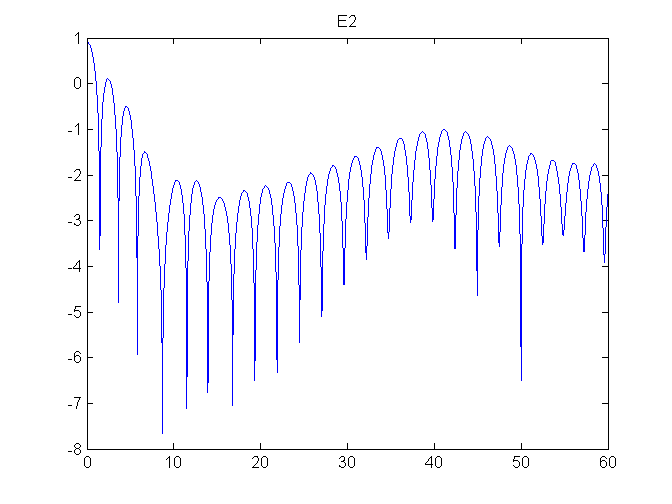},
\includegraphics[width=3in,clip]{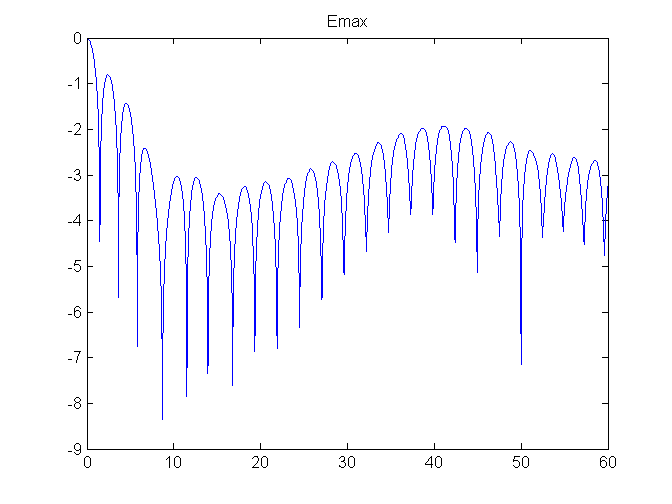} \\
\includegraphics[width=3in,clip]{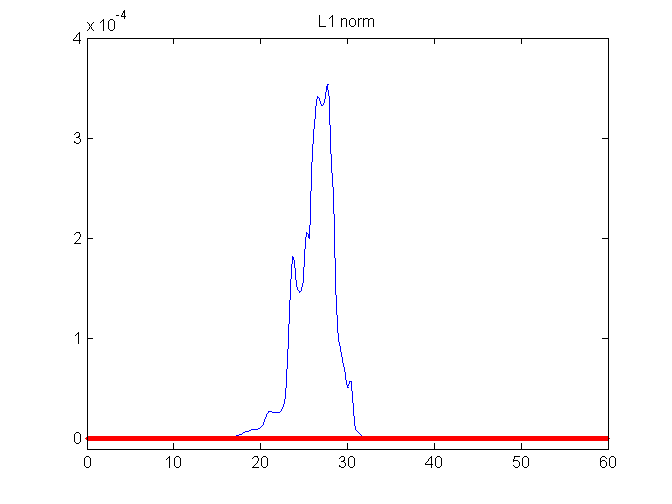},
\includegraphics[width=3in,clip]{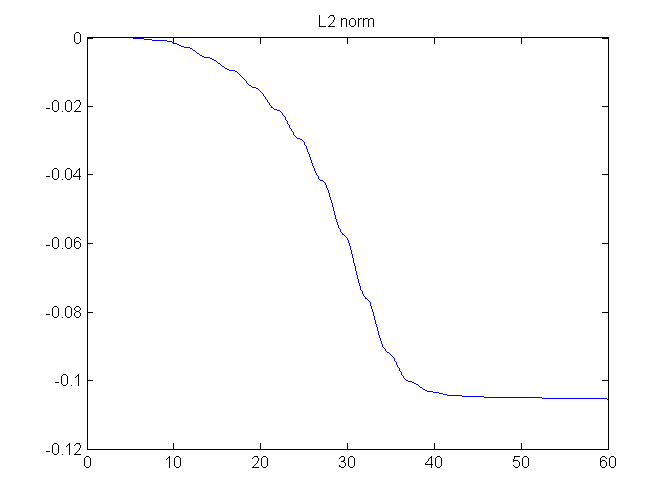}\\
\includegraphics[width=3in,clip]{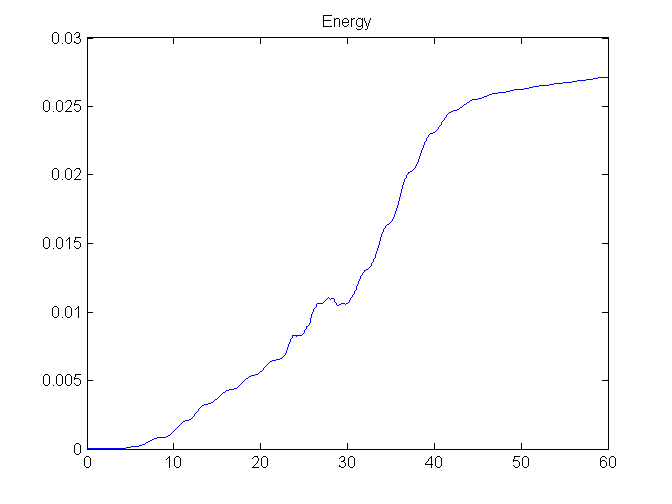},
\includegraphics[width=3in,clip]{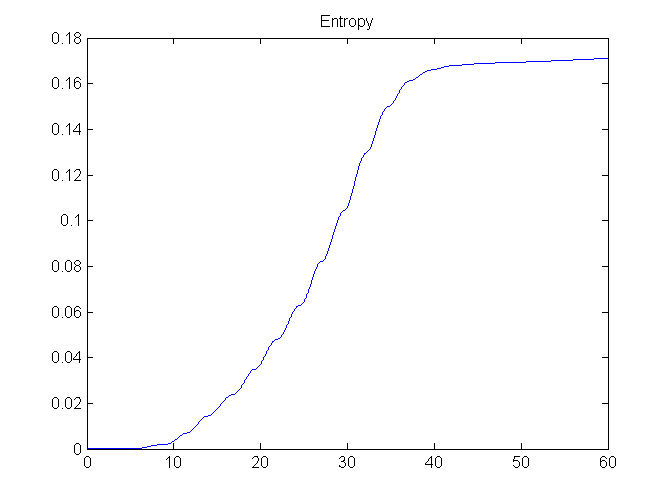}
\caption{Strong Landau damping. Time evolution of the electric field in $L^2$ norm and $L^\infty$ norm (top), discrete $L^1$ norm and $L^2$ norm (middle), kinetic
energy and entropy (bottom). Mesh: $128 \times 128$. The straight red line indicates mass conservation.}
\label{fig22}
\end{figure}
\end{exa}

\begin{exa} (Two stream instability)
Now we consider the two stream instability problem, with an unstable initial distribution
function given by:
\beq
\label{2stream1}
f(t=0,x,v)=\f{2}{7\sqrt{2\pi}}(1+5v^2)(1+\alpha((\cos(2kx)+\cos(3kx))/1.2+\cos(kx))\exp(-\f{v^2}{2})
\eeq
where $\alpha=0.01$ and $k=0.5$. We plot the numerical solution at $T=53$ in Fig. \ref{fig23}. While the time evolution of electric field in $L^2$ norm and $L^\infty$ norm,
the relative derivation of the discrete $L^1$ norm, $L^2$ norm, kinetic energy and
entropy are shown in Fig. \ref{fig24}.

\begin{figure}
\centering
\includegraphics[width=3in,clip]{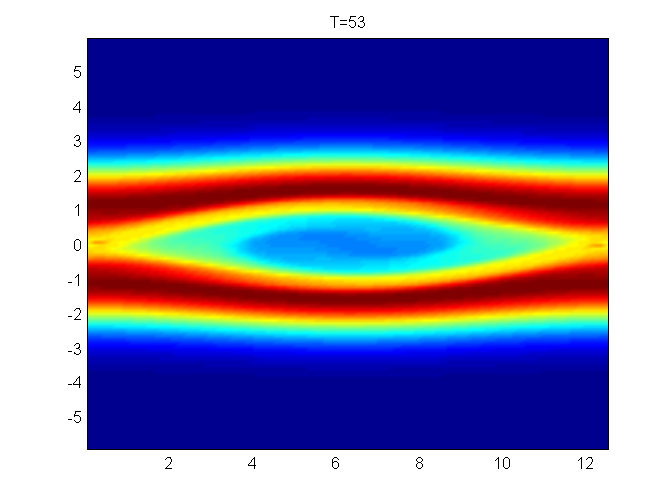},
\includegraphics[width=3in,clip]{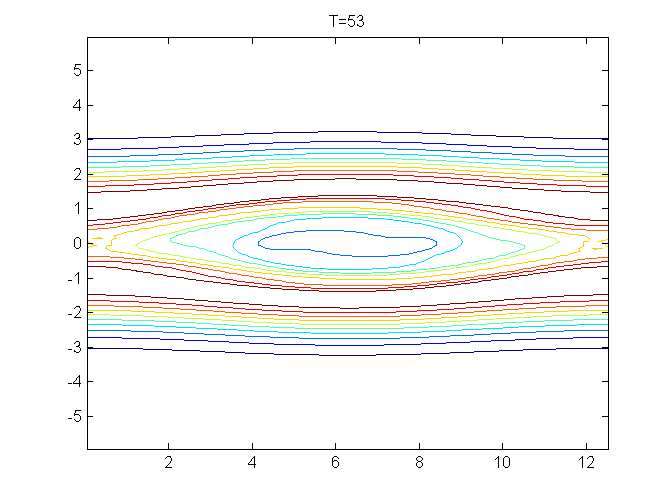}
\caption{Two stream instability at $T=53$. Mesh: $128 \times 128$. Contour plots: 10 equally spaced lines.}
\label{fig23}
\end{figure}

\begin{figure}
\centering
\includegraphics[width=3in,clip]{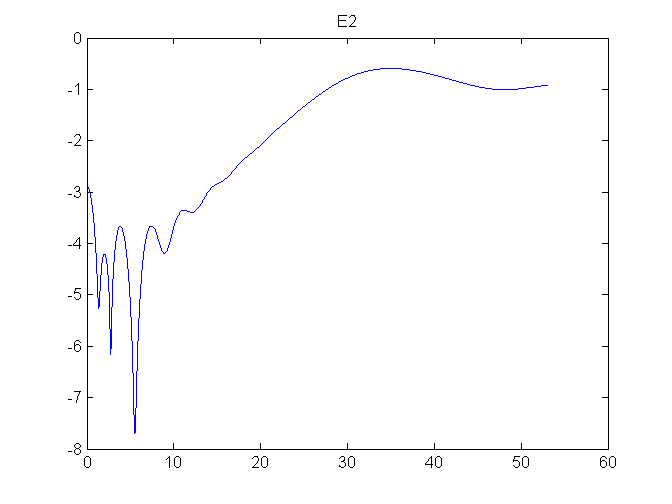},
\includegraphics[width=3in,clip]{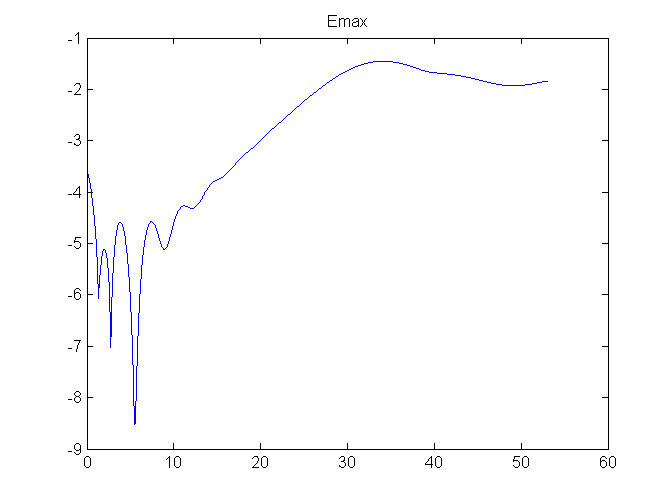} \\
\includegraphics[width=3in,clip]{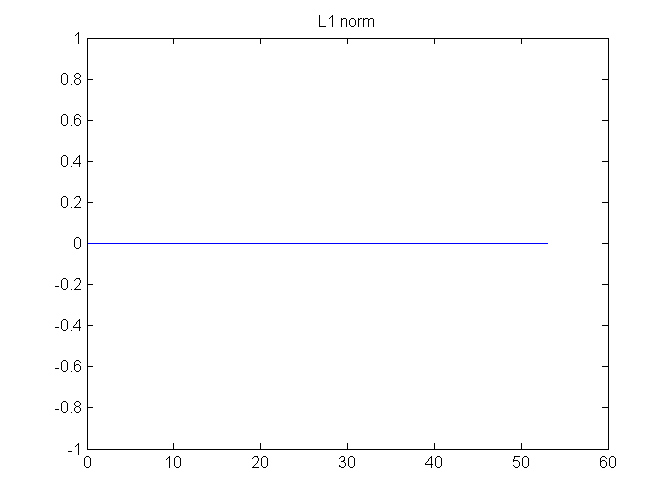},
\includegraphics[width=3in,clip]{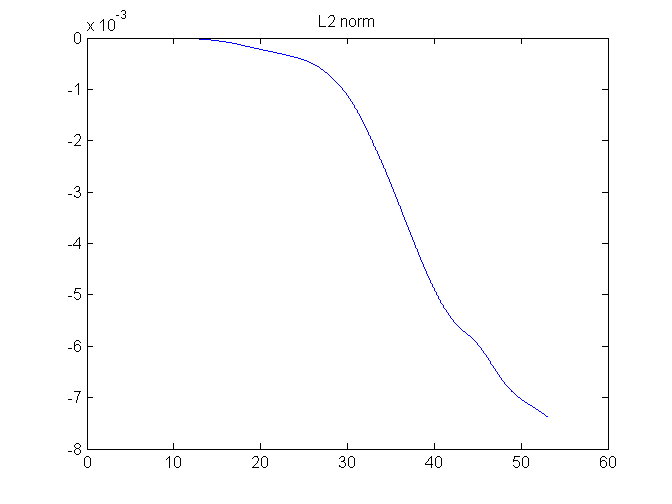}\\
\includegraphics[width=3in,clip]{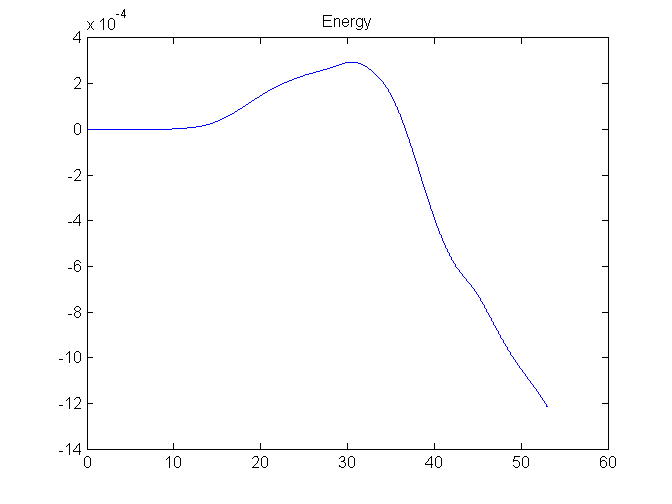},
\includegraphics[width=3in,clip]{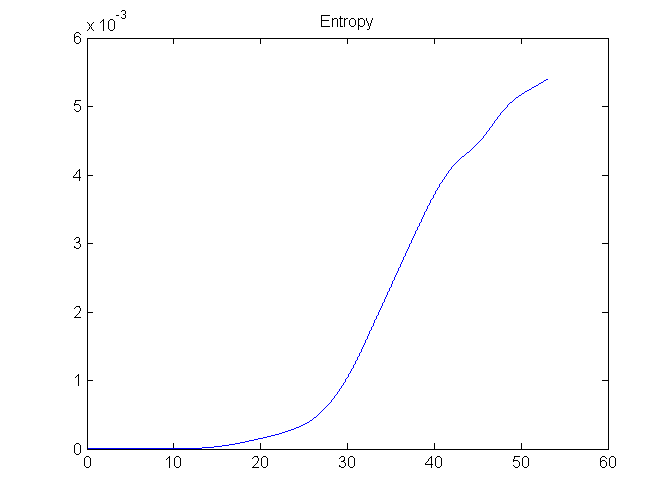}
\caption{Two stream instability. Time evolution of the electric field in $L^2$ norm and $L^\infty$ norm (top), discrete $L^1$ norm and $L^2$ norm (middle), kinetic
energy and entropy (bottom). Mesh: $128 \times 128$.}
\label{fig24}
\end{figure}
\end{exa}

\begin{exa} (Symmetric two stream instability)
We consider the symmetric two stream instability with the initial condition:
\beq
\label{2stream2}
f(t=0,x,v)=\f{1}{2v_{th}\sqrt{2\pi}}\left[\exp\left(-\f{(v-u)^2}{2v_{th}^2}\right)
+\exp\left(-\f{(v+u)^2}{2v_{th}^2}\right)\right](1+\alpha\cos(kx))
\eeq
with $\alpha=0.05$, $u=0.99$, $v_{th}=0.3$ and $k=\f{2}{13}$.
We plot the numerical solution at $T=70$ in Fig. \ref{fig25}. The time evolution
of the electric field in $L^2$ norm and $L^\infty$ norm, the relative derivation of the discrete $L^1$ norm, $L^2$ norm, kinetic energy and
entropy are reported in Fig. \ref{fig26}. Similarly, the mass conservation is indicated by the straight line on the bottom
in the $L^1$ norm figure.

\begin{figure}
\centering
\includegraphics[width=3in,clip]{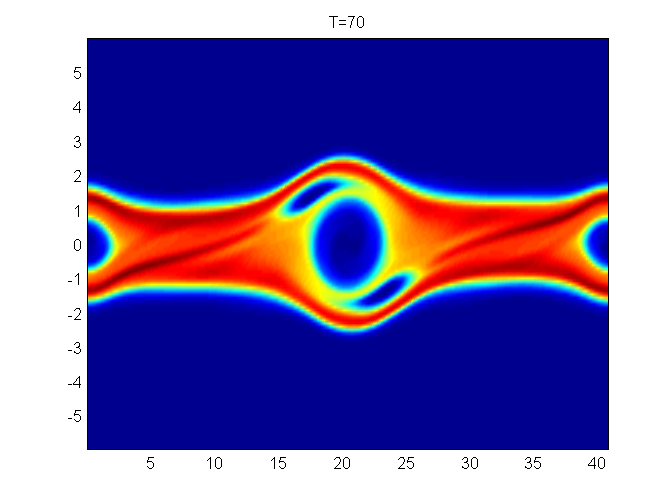},
\includegraphics[width=3in,clip]{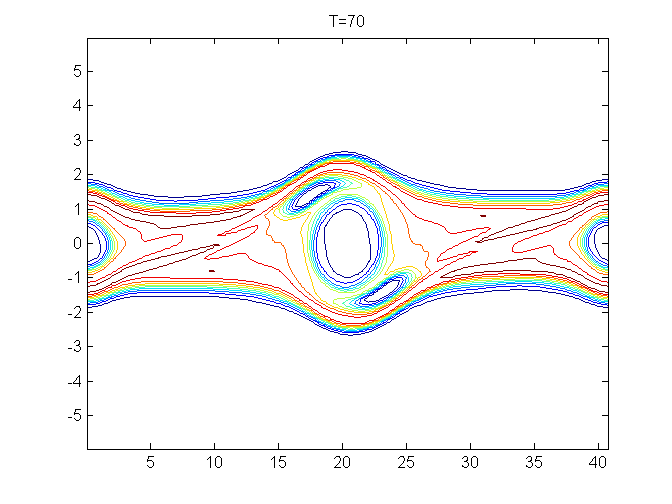}
\caption{Symmetric two stream instability at $T=70$. Mesh: $256 \times 128$. Contour plots: 10 equally spaced lines.}
\label{fig25}
\end{figure}

\begin{figure}
\centering
\includegraphics[width=3in,clip]{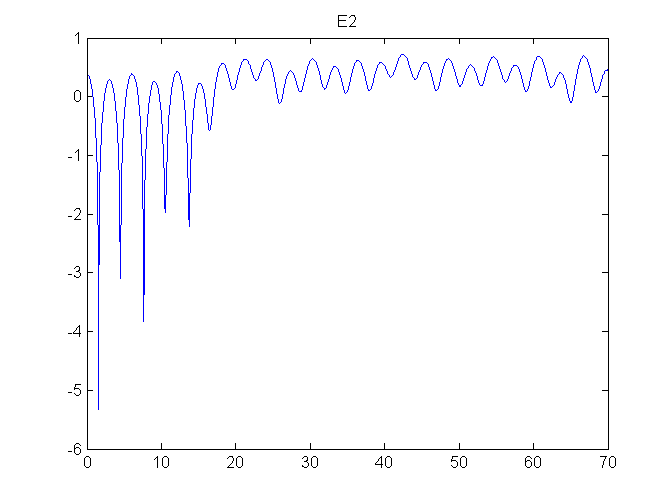},
\includegraphics[width=3in,clip]{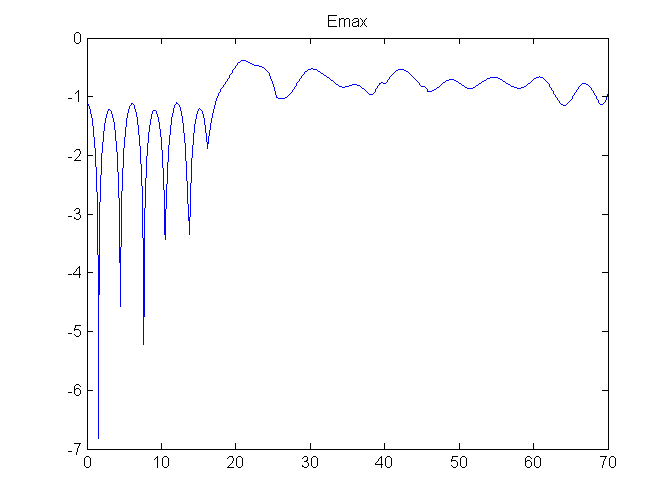} \\
\includegraphics[width=3in,clip]{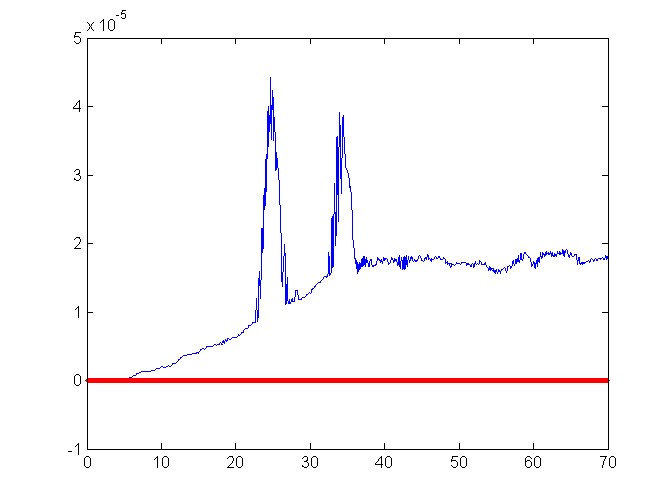},
\includegraphics[width=3in,clip]{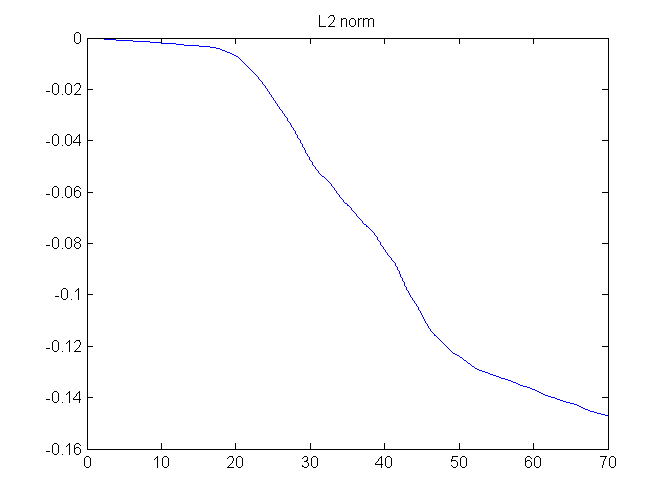}\\
\includegraphics[width=3in,clip]{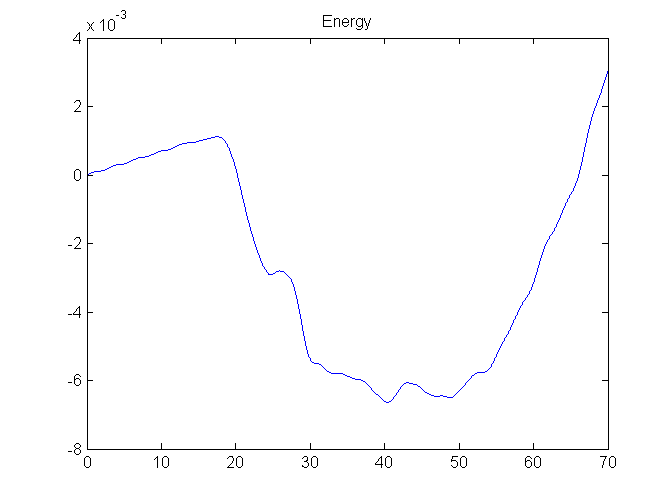},
\includegraphics[width=3in,clip]{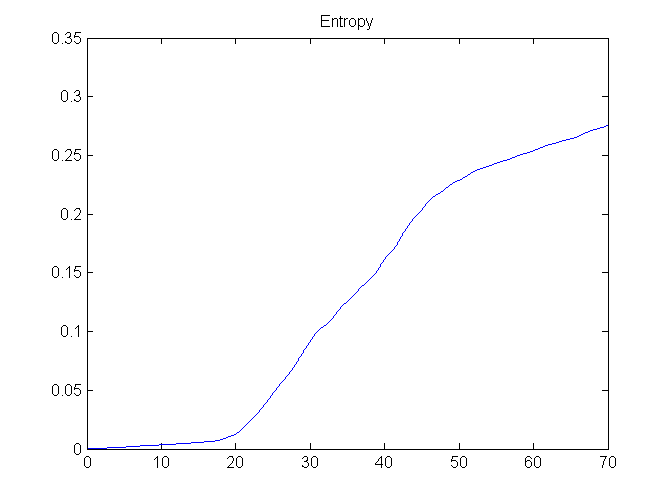}
\caption{Symmetric two stream instability. Time evolution of the electric field in $L^2$ norm and $L^\infty$ norm (top), discrete $L^1$ norm and $L^2$ norm (middle), kinetic
energy and entropy (bottom). Mesh: $256 \times 128$. The straight red line indicates mass conservation.}
\label{fig26}
\end{figure}
\end{exa}

\subsection{The guiding center Vlasov model}

Consider the guiding center approximation of the 2D Vlasov model \cite{yang2014conservative, frenod2015long},
\beq
 \frac{\partial \rho}{\partial t} + E_2  \frac{\partial \rho}{\partial x} - E_1 \frac{\partial \rho}{\partial y} =0,
 \eeq
 or equivalently in a conservative form as
 \beq
 \frac{\partial}{\partial t} \rho+ \frac{\partial }{\partial x} (\rho E_2)+ \frac{\partial }{\partial y}(-\rho E_1) =0,
 \eeq
 where ${\bf E} = (E_1, E_2) = -\nabla \Phi$ with $\Phi$ determined from the Poisson's equation
 $$\bigtriangleup \Phi = -\rho.$$
We assume a uniform set of 2D grid points as specified in eq. \eqref{eq: 2dgrid}. The equations for tracking characteristics emanating from a grid point $(x_i, y_j)$ at some future time $t^{n+1}$ (without loss of generality),
\beq
\label{eq: guiding_char}
\frac{dx(t)}{dt} = E_2, \quad \frac{dy(t)}{dt} =  -E_1, \quad x(t^{n+1}) = x_i, \quad y(t^{n+1}) = y_j.
\eeq
Below we generalize the characteristics tracing procedures in \cite{qiu_russo_2016} to the guiding center model,
which can be directly applied to the incompressible Euler equations in the following subsection. In particular for the system \eqref{eq: guiding_char}, we propose a scheme to locate the foot of characteristics $(x^{\star}_{i,j}, y^{\star}_{i,j})$ at $t^n$. Once the foot of characteristic is located, then a 2D interpolation procedure can be employed to approximate the solution value $\rho(x^{\star}_{i,j}, y^{\star}_{i,j}, t^n)$. We remark that solving \eqref{eq: guiding_char} with high order temporal accuracy is challenging. Especially,
the ${\bf E}$ depends on the unknown function $\rho$ via the 2-D Poisson's equation in a global rather than a local fashion, and it is
difficult to evaluate ${\bf E}$ for some intermedia time stages, i.e. Runge-Kutta methods cannot be used directly.

In our notations, the superscript $^n$ denotes the time level, the subscripts $i$ and $j$ denote the location at $(x_i, y_j)$. e.g. $E^n_{1, i, j} = E_1(x_i, y_j, t^n)$. The superscript $^{(p)}$ denotes the formal order of {\em temporal} approximation. For example, in eq. \eqref{eq: x_v_1} below, $x^{n, (1)}_{i,j}$ (or $y^{n, (1)}_{i,j}$) approximates $x_{i,j}^\star$ (or $y_{i,j}^\star$) with first order. $\frac{d}{dt} = \frac{\partial}{\partial t} +  \frac{\partial x}{\partial t}\frac{\partial}{\partial x}+ \frac{\partial y}{\partial t}\frac{\partial}{\partial y}$
denotes the material derivative along characteristics. We use a spectrally accurate fast Fourier transform (FFT) for solving the 2-D Poisson's equation \eqref{eq: poisson}.

We start from a first order scheme for tracing characteristics \eqref{eq: guiding_char}, by letting
\beq
\label{eq: x_v_1}
x^{n, (1)}_{i, j} = x_i - E_2(x_i, y_j, t^n) \Delta t; \quad y^{n, (1)}_{i, j} = y_j + E_1(x_i, y_j, t^n)  \Delta t.
\eeq
They are first order approximations to $x_{i, j}^\star$ and $y_{i, j}^\star$.
Let
\beq
\rho^{n+1, (1)}_{i, j} = \rho(x^{n, (1)}_{i, j}, y^{n, (1)}_{i, j}, t^n),
\eeq
which can be obtained by a high order spatial interpolation. Based on $\{\rho^{n+1, (1)}_{i, j}\}$, we can compute
\[
{\bf E}^{n+1, (1)}_{i, j} = (E^{n+1, (1)}_{1, i, j}, E^{n+1, (1)}_{2, i, j}),
\]
by using FFT based on the 2-D Poisson's equation \eqref{eq: poisson}.
Note that  ${\bf E}^{n+1, (1)}_{i, j}$ approximates ${\bf E}^{n+1}_{i, j}$ with first order temporal accuracy.

A second order scheme can be built upon the first order one, by letting
 \begin{align}
 \label{eq: x_2}
 & x^{n, (2)}_{i, j} = x_i - \frac12\left (E^{n+1, (1)}_{2, i, j} + E_2(x^{n, (1)}_{i, j}, y^{n, (1)}_{i, j}, t^n)\right) \Delta t,
  \\
  \label{eq: v_2}
  &y^{n, (2)}_{i, j} = y_j + \frac12\left (E^{n+1, (1)}_{1, i, j} + E_1(x^{n, (1)}_{i, j}, y^{n, (1)}_{i, j}, t^n)\right) \Delta t.
 \end{align}
 Here ${\bf E}(x^{n, (1)}_{i, j}, y^{n, (1)}_{i, j}, t^n)$ can be approximated by a high order spatial interpolation.
$(x^{n, (2)}_{i, j}, y^{n, (2)}_{i, j})$ can be shown to be second order approximations to $(x_{i, j}^\star, y_{i, j}^\star)$ by a local truncation error analysis.

Finally, a third order scheme can be designed based on a second order one, by letting
\begin{align}
\label{eq: x_3}
  x^{n, (3)}_{i, j} = x_i - E^{n+1, (2)}_{2, i, j} \Delta t +& \frac{\Delta t^2}2 \Big(
  \frac23
  (\frac{d E_2}{dt})^{n+1, (2)}_{i, j}+ \frac13
  \frac{dE_2}{dt}(x_{i, j}^{n, (2)}, y_{i, j}^{n, (2)},t^n)
   \Big);  \\
  \label{eq: v_3}
  y^{n, (3)}_{i, j} = y_j + E^{n+1, (2)}_{1, i, j} \Delta t - &\frac{\Delta t^2}2 \Big(
  \frac23
  (\frac{d E_1}{dt})^{n+1, (2)}_{i, j}+ \frac13 \frac{d E_1}{dt}(x_{i, j}^{n, (2)}, y_{i, j}^{n, (2)},t^n)
   \Big);
 \end{align}
which are third order approximations to $x_{i, j}^\star$ and $y_{i, j}^\star$, see Proposition~\ref{prop: order3} below.
Here
\beq
\label{eq: material_der}
\frac{d}{dt} E_s = \frac{\partial E_s}{\partial t}+\frac{\partial E_s}{\partial x} E_2 - \frac{\partial E_s}{\partial y} E_1, \quad s=1, 2
\eeq
are material derivatives along characteristics. Notice that on the r.h.s. of eq. \eqref{eq: material_der}, the partial derivatives are not explicitly given. The spatial derivative terms can be approximated by high order spatial approximations, while the time derivative term $\frac{\partial {\bf E}}{\partial t}$ can be approximated by utilizing the Vlasov equation. In particular, taking partial time derivative of the
2-D Poisson's equation gives
\begin{equation}
\Delta \phi_t=-(E_2\rho)_x+(E_1\rho)_y.
\label{eq: poisson2}
\end{equation}
After obtaining ${\bf E}$ by solving the original Poisson's equation \eqref{eq: poisson}, the right hand side of \eqref{eq: poisson2} can be constructed by a high order central finite difference scheme, e.g., 6th
order central finite difference scheme. Then we can solve \eqref{eq: poisson2} by FFT to get $\frac{\partial {\bf E}}{\partial t}= -(( \phi_t )_x, ( \phi_t )_y)$.
With such a procedure, both $\frac{\partial {\bf E}}{\partial t}(x_{i, j}^{n, (2)}, y_{i, j}^{n, (2)},t^n) $ and  $\left(\frac{\partial {\bf E}}{\partial t}\right)^{n+1,(2)}_{i,j}$ can be obtained.

\begin{prop}
\label{prop: order3}
$x^{n, (3)}_{i, j}$ and $y^{n, (3)}_{i, j}$ constructed in equations \eqref{eq: x_3}-\eqref{eq: v_3} are third order approximations to $x_{i, j}^\star$ and $y_{i, j}^\star$ in time.
\end{prop}
\noindent
{\em Proof.}
It can be checked by Taylor expansion
\beqa
x_{i, j}^\star &=& x_i - \frac{d x}{dt}(x_i, y_j, {t^{n+1}}) {\Delta t} + \left(\frac23 \frac{d^2 x}{dt^2}(x_i, y_j, {t^{n+1}}) + \frac13 \frac{d^2 x}{dt^2}(x_{i,j}^\star, y_{i,j}^\star, {t^{n}})\right) \frac{\Delta t^2}{2}
+\mathcal{O}(\Delta t^4) \nonumber\\
&=& x_i -  E^{n+1}_{2, i, j} {\Delta t} +
\left(
\frac23 \frac{d E_2}{d t}(x_i, y_j, t^{n+1})
+\frac13 \frac{d E_2}{d t}(x^\star_{i, j}, y^\star_{i, j}, t^{n})
\right)
\frac{\Delta t^2}{2} +
 \mathcal{O}(\Delta t^4) \nonumber\\
&=& x_i -
(E^{n+1, (2)}_{2, i, j} + \mathcal{O}(\Delta t^3)) {\Delta t} +
 \left(
 \frac23  (\frac{dE_2}{d t})^{n+1, (2)}_{i, j}\right.\\ \noindent
 &&\left.+ \frac13 \frac{dE_2}{dt}(x^{n, (2)}_{i, j}, y^{n, (2)}_{i, j}, t^{n})+ \mathcal{O}(\Delta t^3)
 \right) \frac{\Delta t^2}{2}
+\mathcal{O}(\Delta t^4) \nonumber\\
&\stackrel{\eqref{eq: x_3}}{=}& x^{n, (3)}_i +  \mathcal{O}(\Delta t^4). \nonumber
\eeqa
The second last equality is due to the fact that a second order scheme (with superscript $(2)$) gives locally third order approximations.
Hence $x^{n, (3)}_{i, j}$ (similarly $y^{n, (3)}_{i, j}$) is a fourth order approximation to $x_{i, j}^\star$ (similarly $y_{i,j}^\star$) locally in time for one time step. The approximation is third order in time globally.
$\Box$.

\begin{exa} (Kelvin-Helmholtz instability problem). This example is the 2-D guiding center model problem with the initial condition
\begin{equation}
\rho_0(x,y)=\sin(y)+0.015\cos(kx)
\end{equation}
and periodic boundary conditions on the domain $[0,4\pi]\times[0,2\pi]$. We let $k=0.5$, which
will create a Kelvin-Helmholtz instability.

For this example, we show the surface and contour plots for the solution at $T=40$ in Fig. \ref{fig33}, similar to the results
in \cite{frenod2015long}. The mesh size is $128\times128$.

\begin{figure}
\begin{center}
\includegraphics[height=2.5in,width=3.0in]{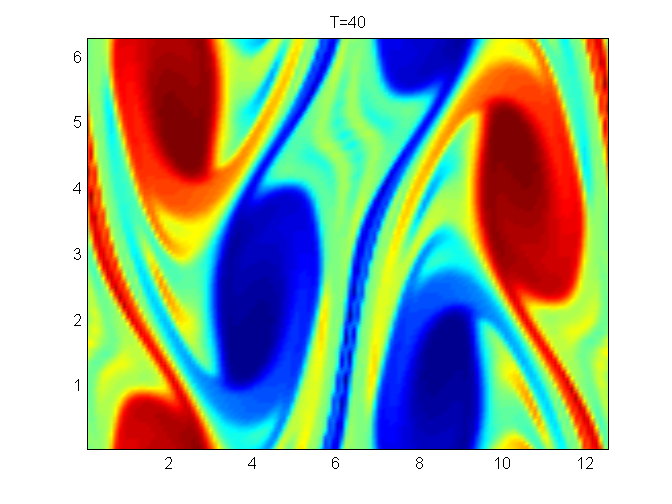}
\includegraphics[height=2.5in,width=3.0in]{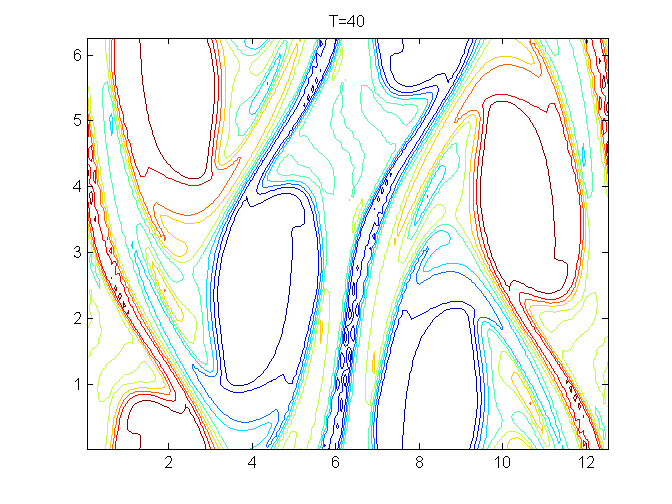}
\end{center}
\caption{Kelvin-Helmholtz instability problem. Mesh size $128\times128$. $T=40$. Contour plots: 10 equally spaced lines. }
\label{fig33}.
\end{figure}


\end{exa}

\subsection{Incompressible Euler equation}

\begin{exa}
\label{ex67}
We first consider the incompressible Euler system on the domain $[0, 2\pi]\times[0,2\pi]$ with an initial condition $\omega_0(x,y)=-2\sin(x)\sin(y)$.
The exact solution will stay stationary with $\omega(x,y,t)=-2\sin(x)\sin(y)$.
Similarly as in Table \ref{tab1} and Table \ref{tab1t}, the $5$th order spatial accuracy and
3rd order temporal accuracy are clearly observed in Table \ref{tab3} and Table \ref{tab3t} respectively.
Here for the temporal accuracy, 7th order linear interpolation and linear reconstruction are used.

\begin{table}[h]
\begin{center}
\caption{{Errors and orders for the incompressible Euler equation in Example \ref{ex67}. $T=1.2$.}
}
\bigskip
\begin{tabular}{|c | c|c|c| c|c|}
\hline
\cline{1-5} $N_x \times N_y$  &{$20\times 20$} &{$40\times 40$} &{$60\times 60$}&{$80\times 80$}&{$100\times100$}  \\
\hline
\cline{1-5}  $L^1$ error &1.00E-2&3.01E-4 &3.80E-5&8.79E-6&2.84E-6\\
\hline
 order &--&5.06&5.10&5.09&5.07\\
\hline
\end{tabular}
\label{tab3}
\bigskip
\caption{{Errors and orders for the incompressible Euler equation in Example \ref{ex67}. $N_x=N_y=128$. $T=1$.}}
\bigskip
\begin{tabular}{|c | c|c|c| c|c|}
\hline
\cline{1-4} $CFL$  &{$1.15$} &{$1.05$} &{$0.95$}&{$0.85$} \\
\hline
\cline{1-4}  $L^1$ error &4.80E-9&3.65E-9 &2.69E-9&1.92E-9\\
\hline
 order &--&3.03&3.04&3.04\\
\hline
\end{tabular}
\label{tab3t}
\end{center}
\end{table}


\end{exa}

\begin{exa} (The vortex patch problem).
In this example, we consider the incompressible Euler equations with the
initial condition given by
\begin{equation}
\omega_0(x,y)=
\begin{cases}
-1, \qquad & \frac{\pi}{2} \le x \le \frac{\pi}{4}\le y \le \frac{3\pi}{4}; \\
1, \qquad & \frac{\pi}{2} \le x \le \frac{5\pi}{4} \le y \le \frac{7\pi}{4}; \\
0, \qquad & \text{otherwise}.
\end{cases}
\end{equation}
We show the surface and contour plots of $\omega$ at $T=5$ in Fig. \ref{fig31}. The mesh size is $128\times128$.

\begin{figure}
\begin{center}
\includegraphics[height=2.5in,width=3.0in]{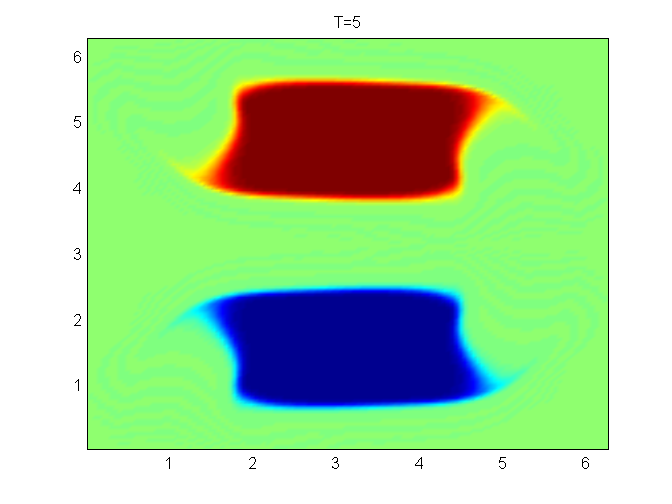}
\includegraphics[height=2.5in,width=3.0in]{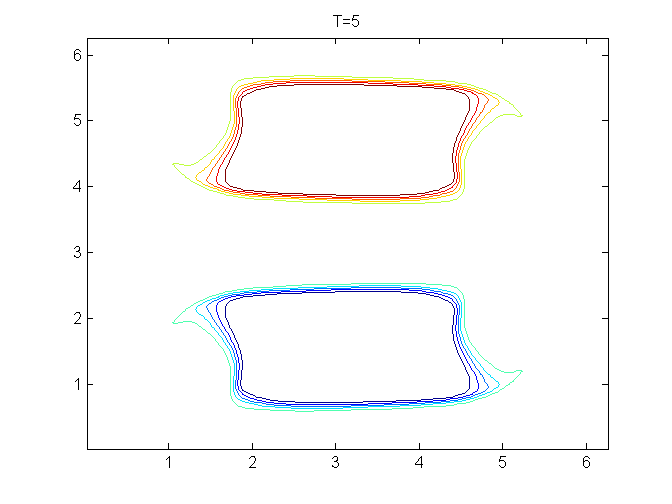}
\end{center}
\caption{Vortex patch problem. Mesh size $128\times128$. $T=5$. Contour plot: 10 equally spaced lines.}
\label{fig31}
\end{figure}

\end{exa}

\begin{exa} (Shear flow problem).
This example is the same as above but with following
initial conditions
\begin{equation}
\omega_0(x,y)=
\begin{cases}
\delta \cos(x)-\frac{1}{\rho}sech^2((y-\pi/2)/\rho)^2, \qquad &y\le \pi; \\
\delta \cos(x)+\frac{1}{\rho}sech^2((3\pi/2-y)/\rho)^2,\qquad & y>\pi.
\end{cases}
\end{equation}
where $\delta=0.05$ and $\rho=\frac{\pi}{15}$.
We show the surface and contour plots of $\omega$ at $T=6$ (top) and $T=8$ (bottom) in Fig. \ref{fig32}. The mesh size is $128\times128$.

\begin{figure}
\begin{center}
\includegraphics[height=2.5in,width=3.0in]{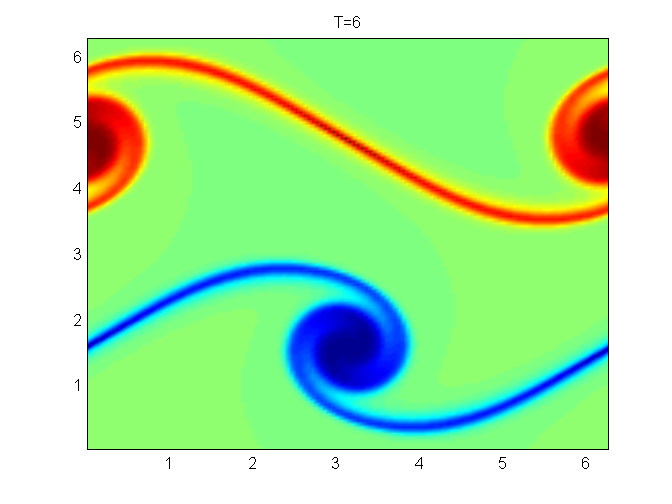}
\includegraphics[height=2.5in,width=3.0in]{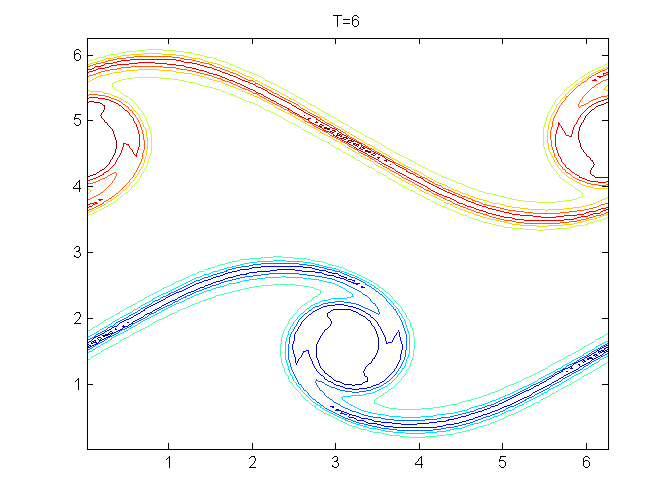}\\
\includegraphics[height=2.5in,width=3.0in]{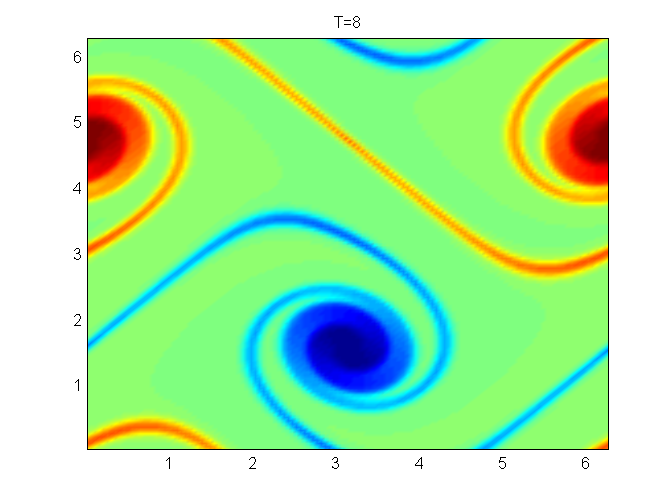}
\includegraphics[height=2.5in,width=3.0in]{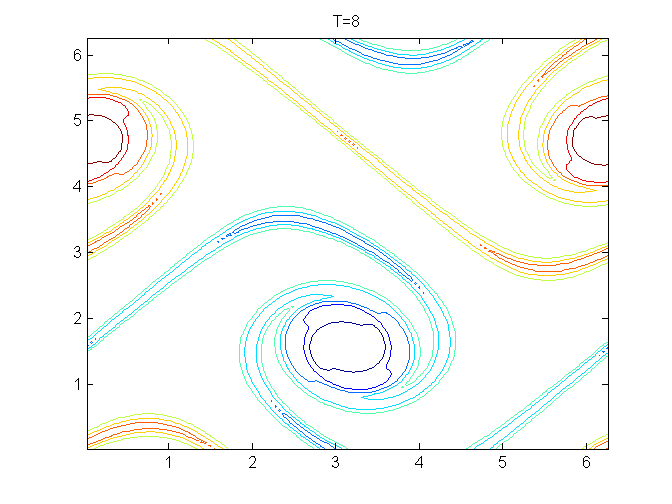}
\end{center}
\caption{Shear flow problem. Mesh size $128\times128$. $T=6$ (top) and $T=8$ (bottom). Contour plots: 10 equally spaced lines.}
\label{fig32}
\end{figure}

\end{exa}

%% file: conclusion.tex
\section{Conclusion}
\label{sec7}
\setcounter{equation}{0}
\setcounter{figure}{0}
\setcounter{table}{0}

In this paper, we propose a conservative semi-Lagrangian finite difference scheme based on a flux difference formulation. We investigate its numerical stability from the linear ODE and PDE point of view via Fourier analysis. The upper bound of time step constraints have been found in the linear setting and have been numerically verified. These upper bounds are only slightly greater than those from the Eulerian approach, unfortunately. The schemes are applied to passive transport problems as well as nonlinear Vlasov systems and the incompressible Euler system to showcase its effectiveness. A new characteristics tracing procedure for the guiding center Vlasov system and incompressible Euler system is proposed, mimicking the characteristic tracing mechanism in \cite{qiu_russo_2016}.